\documentclass[11pt,a4paper]{article}

\usepackage{amsmath,amssymb}
\usepackage{color}
\newcommand{\R}{\mathbb{R}}
\newcommand{\N}{\mathbb{N}}
\newcommand{\Z}{\mathbb{Z}}
\newcommand{\Fs}{\mathbb{F}^*}

\newcommand{\Cs}{\mathbb{C}^*}
\newcommand{\Ds}{\mathbb{D}^*}
\newcommand{\cu}{\underline{c}}
\newcommand{\co}{\overline{c}}
\def\cA{{\mathcal A}}
\def\cC{{\mathcal C}}

\def\cS{{\mathcal S}}
\newcommand{\ee}{\varepsilon}
\renewcommand{\aa}{\alpha}
\renewcommand{\div}{{\rm div}\,}

\newcommand{\Frac}{\displaystyle \frac}

\newcommand{\Sum}{\displaystyle \sum}

\def\d{\partial}
\def\ddl{\dot \Delta_l}
\def\ddj{\dot \Delta_j}
\def\ddq{\dot \Delta_q}

\def\tilde{\widetilde}
\def\hat{\widehat}

\newcommand{\n}{\nabla}
\newcommand{\fd}{\frac{d}{2}}
\newcommand{\fdp}{\frac{d}{p}}
\newcommand{\p}{\partial}

\newcommand{\qe}{q_\ee}
\newcommand{\ue}{u_\ee}
\newcommand{\qa}{q_\aa}
\newcommand{\ua}{u_\aa}
\newcommand{\ca}{c_\aa}
\newcommand{\phie}{\phi_\ee}
\newcommand{\psia}{\psi_\aa}
\newcommand{\dq}{\delta q}
\newcommand{\du}{\delta u}

\newtheorem{thm}{Theorem}
\newtheorem{lem}{Lemma}

\newtheorem{prop}{Proposition}
\newtheorem{defi}{Definition}

\newtheorem{rem}{Remark}

\title{Local in time results for local and non-local capillary Navier-Stokes systems with large data}

\author{Fr\'ed\'eric Charve\footnote{Universit\'e Paris-Est Cr\'eteil, Laboratoire d'Analyse et de Math\'ematiques Appliqu\'ees (UMR 8050), 61 Avenue du G\'en\'eral de Gaulle, 94 010 Cr\'eteil Cedex (France). E-mail: frederic.charve@u-pec.fr}}

\date{}
\begin{document}

\maketitle

\begin{abstract} In this article we study three capillary compressible models (the classical local Navier-Stokes-Korteweg system and two non-local models) for large initial data, bounded away from zero, and with a reference pressure state $\bar{\rho}$ which is not necessarily stable ($P'(\bar{\rho})$ can be non-positive). 
We prove that these systems have a unique local in time solution and we study the convergence rate of the solutions of the non-local models towards the local Korteweg model. The results are given for constant viscous coefficients and we explain how to extend them for density dependant coefficients.
\end{abstract}
\section{Introduction}
\subsection{Presentation of the systems}
In this article we are interested in the dynamics of a liquid-vapor mixture in the setting of the Diffuse Interface (DI) approach: between the two phases lies a thin region of continuous transition and the phase changes are read through the variations of the density (with for example a Van der Waals pressure). Unfortunately the basic models provide an infinite number of solutions and in order to select the physically relevant solutions, following Van der Waals and Korteweg, one penalizes the high variations of the density thanks to capillary terms related to surface tension.

We first consider the classical compressible Navier-Stokes system (NSK) endowed with an internal local capillarity (this is why we call (NSK) the local Korteweg system). This system, first considered by Korteweg and renewed by Dunn and Serrin, reads:
$$
\begin{cases}
\begin{aligned}
&\d_t\rho+\div (\rho u)=0,\\
&\d_t (\rho u)+\div (\rho u\otimes u)-\cA u+\nabla(P(\rho))=\kappa\rho\nabla D[\rho],\\
\end{aligned}
\end{cases}
$$
where $\rho$ and $u$ denote the density and the velocity ($\rho$ is a non-negative function and $u$ is a vector-valued function defined on $\R_+\times\R^d$). The general diffusion operator is defined as follows:
$$
\cA u = \div (2\mu(\rho) Du) +\nabla (\lambda(\rho) \div u),
$$
where $2Du =^t\nabla u +\nabla u$. For simplicity we will present the results in the case of constant viscosity coefficients (we refer to the end of the article for general coefficients) so that
$$
\cA u= \mu \Delta u+ (\lambda+\mu)\nabla \div u, \quad \mbox{with} \quad \mu>0 \quad \mbox{and} \quad \nu=\lambda+ 2 \mu >0.
$$
In the classical Korteweg system $(NSK)$, the capillary term is defined by $\div K$, where the general Korteweg tensor is given by:
$$
K(\rho)= \frac{\kappa(\rho)}{2} (\Delta \rho^2 -|\nabla \rho|^2) I_d -\kappa(\rho) \nabla \rho \otimes \nabla \rho,
$$
The coefficient $\kappa$ may depend on $\rho$ but in this article it is chosen constant and then the capillary term turns into $\kappa \rho \nabla D[\rho]$ with (see \cite{3DS}):
$$D[\rho]=\Delta\rho,$$
The solutions of this system are much more regular than those of the classical compressible system (NSC) (that is for $D[\rho]=0$). We refer to \cite{3K, 3DS, DD} for more informations about this model.

Another way of selecting the physical solutions consists in defining a non-local capillary term involving the density through a convolution and only one derivative (compared to the numerical difficulties generated by the previous local capillary term with derivatives of order 3). In the non-local Korteweg system $(NSRW)$, introducing $\phi$, called interaction potential, which satisfies the following conditions
\begin{equation}
 (|.|+|.|^2)\phi(.)\in L^1(\R^d)\mbox{, }\quad\int_{\R^d}\phi(x)dx=1,\quad\phi\mbox{ even, and }\phi\geq0,
\label{potentielcond}
\end{equation}
then $D[\rho]$ is the non-local term given by:
$$D[\rho]=\phi*\rho-\rho.$$

Computing the Fourier transform of these capillary terms, $(\hat{\phi}(\xi)-1) \hat{\rho}(\xi)$ in the non-local model, and $-|\xi|^2 \hat{\rho}(\xi)$ in the local model, a natural question arises which is to study the closedness of the solutions of these models when $\hat{\phi}(\xi)$ is formally "close" to $1-|\xi|^2$. To answer this question, we chosed in \cite{CH} a particular interaction potential and considered the following non-local system: 
$$
\begin{cases}
\begin{aligned}
&\d_t\rho_\ee+\div (\rho_\ee \ue)=0,\\
&\d_t (\rho_\ee \ue)+\div (\rho_\ee u\otimes \ue)-\cA \ue+\nabla(P(\rho_\ee))=\kappa \rho_\ee \nabla D[\qe],\\
\end{aligned}
\end{cases}
\leqno{(NSRW_\ee)}
$$
where
$$\begin{cases}
\phi_{\ee}=\frac{1}{\ee^d} \phi(\frac{x}{\ee}) \quad \mbox{with} \quad \phi(x)=\frac{1}{(2 \pi)^d} e^{-\frac{|x|^2}{4}},\\
D[\qe]=\frac{1}{\ee^2}(\phi_\ee*\rho_\ee-\rho_\ee).
\end{cases}
$$
When $\xi$ is fixed, we have $\hat{\phie}(\xi)=e^{-\ee^2|\xi|^2}$, and when $\ee$ is small, $\Frac{\hat{\phi_\ee}(\xi)-1}{\ee^2}$ is close to $-|\xi|^2$.
We refer to \cite{VW, Rohdehdr, 5CR, 5Ro, Has2, Has5, CH, CHVP} for more details. The solutions of the non-local model have a regularity structure closer to what is obtained for system (NSC) but the numerical difficulties seem comparable to (NSK) due to the convolution operator.

For this reason, C. Rohde introduced in \cite{Rohdeorder} a new model, called the order-parameter model, and inspired by the work of D. Brandon, T. Lin and R. C. Rogers in \cite{BLR}.  In this new system the capillary term $\alpha^2 \nabla (c-\rho)$ involves a new variable $c$ called the "order parameter", which is coupled to the density via the following relation coming from the Euler-Lagrange equation from the variational approach ($\aa$ controls the coupling between $\rho$ and $c$):
$$
\Delta c +\alpha^2 (\rho-c)=0,
$$
so that the new system he considered is the following:
$$
\begin{cases}
\begin{aligned}
&\d_t\rho_\aa+\div (\rho_\aa \ua)=0,\\
&\d_t (\rho_\aa \ua)+\div (\rho_\aa \ua\otimes \ua)-\cA \ua+\nabla(P(\rho_\aa))=\kappa \aa^2 \rho_\aa\nabla(\ca-\rho_\aa),\\
&\Delta \ca +\alpha^2 (\rho_\aa -\ca)=0.
\end{aligned}
\end{cases}
\leqno{(NSOP_\aa)}
$$
As emphasized by C. Rohde, from a numerical point of view this system is much more interesting (only one derivative in the capillary tensor which is local), and the additionnal equation for the order parameter is a simple linear elliptic equation that can be easily solved (and numerically fast). Another important feature of this model is that thanks to the capillary term, the momentum equation can be rewritten with a modified pressure: $\tilde{P}(\rho)= P(\rho) +\kappa \aa^2 \rho^2/2$, so that if $\aa$ is large enough the derivative will be positive.

Moreover as explained in \cite{Corder} introducing the following interaction potential: 
$$\psia=\aa^d \phi(\aa \cdot) \quad \mbox{with} \quad \psi(x)=\frac{C_d}{|x|^{\fd-1}} K_{\fd-1}(|x|),$$
where $K_\nu$ denotes the modified Bessel function of the second kind and index $\nu$, then the system can be rewritten in the shape of the previous non-local capillary model:
$$
\begin{cases}
\begin{aligned}
&\d_t\rho_\aa+\div (\rho_\aa \ua)=0,\\
&\d_t (\rho_\aa \ua)+\div (\rho_\aa u\otimes \ua)-\cA \ua+\nabla(P(\rho_\aa))=\rho_\aa \kappa \aa^2 \nabla(\psia*\rho_\aa-\rho_\aa),\\
\end{aligned}
\end{cases}
\leqno{(NSOP_\aa)}
$$
\begin{rem}
\sl{From the previous computations, we immediately get that
$$
\ca-\rho_\aa=(-\Delta+ \aa^2 I_d)^{-1}\Delta \rho_\aa=\psia * \rho_\aa - \rho_\aa,
$$
that is $\ca=\psia * \rho_\aa$. This is why the only choice for the initial order parameter is $c_0=\psia *\rho_0$.
\label{orderparamexpr}
}
\end{rem}
We refer to \cite{BLR, Rohdeorder, Corder} for more details.

We focus here on strong solutions with initial data in critical spaces. Let us recall that a critical space is a space whose corresponding norm has the same scaling invariance as the (NSC) system: if $(\rho(t,x), u(t,x))$ is a solution corresponding to the initial data $(\rho_0(x), u_0(x))$, then for each $\lambda>0$, $(\rho(\lambda^2 t, \lambda x), \lambda u(\lambda^2 t, \lambda x))$ is also a solution, corresponding to the dilated initial data $(\rho_0(\lambda x), \lambda u_0(\lambda x))$, provided that the pressure $P$ has been changed into $\lambda^2 P$. For example the Sobolev space $\dot{H}^\fd (\R^d)$, or the Besov space $\dot{B}_{p,1}^\fdp$ are critical. We refer to \cite{Dinv, Dbook, Dheatglobal, Dheatlocal, CD} for more details.

A natural first step in the study of system (NSC) is to consider initial data in critical spaces close to an equilibrium state $(\overline{\rho},0)$ with $P'(\overline{\rho})>0$. Assuming that the initial density fluctuation $q_0$, defined by $\rho_0= \overline{\rho}(1+q_0)$, is small, we perform the classical change of function $\rho= \overline{\rho}(1+q)$ and expect that $q$ is also small. For simplicity we take $\overline{\rho}=1$, then $q=\rho-1$ is expected to be small, and $\rho$ will be bounded away from zero. Simplifying by $\rho$, the previous systems become (we chosed to drop the subscripts):
$$
\begin{cases}
\begin{aligned}
&\d_t q+ u\cdot\nabla q+ (1+q)\div u=0,\\
&\d_t u+ u\cdot\nabla u -\frac{1}{1+q}\cA u+\frac{P'(1+q)}{1+q}\cdot\nabla q-\kappa \nabla D[q]=0,\\
\end{aligned}
\end{cases}
$$
where
\begin{equation}
\begin{cases}
D[q]=0\mbox{ for }(NSC),\\
D[q]=\Delta q \mbox{ for }(NSK),\\
D[\qe]=\displaystyle{\frac{\phie *\qe-\qe}{\ee^2}}$ \mbox{ for } $(NSRW_\ee),\\
D[\qa]=\aa^2 (\psia * \qa-\qa)$ \mbox{ for } $(NSOP_\aa).
\end{cases}
\label{Capill}
\end{equation}
If we also assume that the equilibrium $(1,0)$ is stable (that is $P'(1)>0$) then it is useful to rewrite the system into:
$$
\begin{cases}
\begin{aligned}
&\d_t q+ u\cdot\nabla q+ (1+q)\div u=0,\\
&\d_t u+ u\cdot\nabla u -\cA u+P'(1)\nabla q-\kappa \nabla D[q]=K(q)\nabla q- I(q) \cA u,\\
\end{aligned}
\end{cases}
$$
where $K$ and $I$ are the following real-valued functions defined on $\R$:
$$
K(q)=\left(P'(1)-\frac{P'(1+q)}{1+q}\right) \quad \mbox{and} \quad I(q)=\frac{q}{q+1}.
$$
If the density fluctuation $q$ is small, then $K(q)$ and $I(q)$ are expected to be small so that the contribution of the right-hand side term should also be small.

Indeed when $q_0$ is small in $\dot{B}_{2,1}^{\fd-1} \cap \dot{B}_{2,1}^\fd$ and $P'(1)>0$, R. Danchin obtained in \cite{Dinv, Dheatglobal} global existence of a solution for the compressible Navier-Stokes system $(NSC)$. The local capillary model has been treated by R. Danchin and B. Desjardins in \cite{DD} and we refer to \cite{Has2, CH, CHVP, Corder} for the same results in the non-local capillary case. In addition we proved that the solutions of these non-local models converge towards the solution of the local model with the same small initial data $q_0$, and provided an explicit rate of convergence. In \cite{CHVP, Corder} we give refined estimates.

When the constant state $(1,0)$ is not assumed to be stable anymore, the best we can obtain are local in time existence results. Such results were first obtained by R. Danchin when $q_0$ is small in $\dot{B}_{2,1}^\fd$ (then again $1+q_0$ is bounded away from zero: no vacuum). We refer to \cite{Dheatlocal, Dbook} for the $(NSC)$ system and to \cite{DD} for the local capillary model. We emphasize that no assumption on the pressure law or on the stability are needed, so that the case of a Van Der Waals pressure law is covered. From \cite{CH, Corder} these existence results can be very easily adapted to the non-local capillary models as well as the convergence, when $\ee$ is small or $\aa$ is large for a finite lifespan $T$.

A much more difficult question is to study these systems whithout any stability assumptions on the state $(1,0)$ and any smallness condition on the initial density fluctuation $q_0$ only assumed to belong to $\dot{B}_{2,1}^\fd$, which is naturally the context when two phases coexist. For example if $\rho_0= (1-\chi) \rho_1 + \chi \rho_2$ with $\rho_{1,2}$ two constants, and $\chi$ is a smooth cut-off function.

The only assumption is that $\rho_0=1+q_0$ is bounded away from zero. As explained more in details later, in this unfavorable case we cannot rely on the a priori estimates used in the previous works because of the terms $q \div u$ and $I(q) \cA u$. When $q$ is small these terms are harmless because their  Besov-norms are easily absorbed by the left-hand side. On the contrary, in the present case $q$ has no reason to be small and none of the previous terms, even in a small intervall of time, can be handled like previously and both of them obtruct any use of the previous estimates. The idea introduced by R. Danchin to deal with large density fluctuation is basically to decouple $(q,u)$ and study a slightly modified equation on the velocity, where thanks to a frequency truncation, a big part of the previous problematic term can be included in the linear system and the rest can be made small and absorbable by the left-hand side. Roughly speaking, instead of studying:
$$
\d_t u+ v\cdot\nabla u -\cA u=-\frac{P'(1+q)}{1+q}\nabla q -\frac{q}{1+q} \cA u,
$$
R. Danchin, studied in \cite{Dtruly}:
$$
\d_t u+ u\cdot\nabla u -\dot{S}_m (\frac{1}{1+q})\cA u=-\frac{P'(1+q)}{1+q} \nabla q +\left( (I_d-\dot{S}_m )\frac{1}{1+q}\right) \cA u.
$$
Using refined estimates (see \cite{Dtruly, Dbook}) on the density fluctuation equation, we can fix $m$ large enough (only depending on the initial data), so that $(I_d-\dot{S}_m )\frac{1}{1+q}$ is small and the last term can be absorbed through estimates for the following equation: 
\begin{equation}
\d_t u+ v\cdot\nabla u -b\cA u=F,
\label{estimub}
\end{equation}
where $b$ is a regular function, bounded away from zero, the first part of the right-hand side being small in a small interval of time.

To the best of our knowledge there are very few similar results in the capillary case or for density dependant coefficients (viscosity, capillarity). For special choices on the viscosity and capillary coefficients, the previous method can be simplified as there is no need for new a priori estimates. For instance we refer to \cite{CHSW} in the case of the shallow water model, that is when $\mu(\rho)=\rho$, $\lambda(\rho)=0$, the system turns into:
$$
\begin{cases}
\begin{aligned}
&\p_{t}q+u\cdot\n q+(1+q){\rm div}u=0,\\
&\p_{t}u+u\cdot\n u-{\mathcal  A}u-2D(u).\n(\ln(1+q))+\n(H(1+q))=0 .
\end{aligned}
\end{cases}
$$
And the problematic term can be decomposed into: $\n\ln(1+q)\cdot D(u)=I+II$, with 
$$  
I=\n\Big(\ln(1+q)-\ln(1+\dot{S}_m q)\Big)\cdot D(u), \ \ \
II=\n\ln(1+\dot{S}_m q)\cdot D(u). 
$$ 
As before, the first one is small for large fixed $m$ (we refer to \cite{Dtruly, Dbook, CHSW}), and the second one will be proved to be small thanks to the smoothness of $\dot{S}_m q$: for $0<\alpha<1$,
\begin{multline}
 \|II\|_{\dot{B}_{2,1}^{\fd-1}} \leq
 \|\n\ln(1+\dot{S}_m q)\|_{\dot{B}_{2,1}^{\fd+\alpha-1}} \|D(u)\|
_{\dot{B}_{2,1}^{\fd-\alpha}}\\
\leq C(\|q\|_{L^\infty}) 2^{m\alpha}
 \|q\|_{\dot{B}_{2,1}^{\fd}} \|u\|_{\dot{B}_{2,2}^{\fd+1-\alpha}},
\end{multline}
which is partly absorbed by the left-hand side when $t\in[0,T]$ with $T$ small enough, the other part being dealt thanks to the Gronwall lemma. We emphasize that this could not have been performed with $q$ instead of $\dot{S}_m q$. Let us precise that in \cite{CHSW} there is an improvement on the assumptions on the initial data velocity which is taken in a bigger space: $u_{0}\in \dot{B}^{\fd-1}_{2,2}\cap \dot{B}^{-1}_{\infty,1}$, $\div u_{0}\in \dot{B}^{\fd-2}_{2,1}$.

We also refer to \cite{Has7, Has8} where B. Haspot obtains local results for large data for the local Korteweg model with special choices of the pressure and the viscosity and capillarity coefficients: $\kappa(\rho)=\frac{1}{\rho}$, $\mu(\rho)=\rho$, $\lambda(\rho)=\rho$ or $0$ and $P(\rho)=\rho$. We emphasize that another important feature of \cite{Has7, Has8} is that the initial data are taken in Besov spaces with third index $2$ or infinite. In these works it is more convenient to introduce $q=\log \rho$ and the effective velocity $v=u+\nabla \ln \rho$ which has important regularity properties. We refer to \cite{Hasarma, CH1} for other use of the notion of effective velocity.

\section{Statement of the results}

The present work is devoted to the study of the capillary models with large data in the case of constant viscosity and capillarity coefficients without any stability assumption on the state $(1,0)$.
Due to the capillary term that either involves too much derivatives in (NSK) or large coefficients in terms of $\ee$ or $\aa$ in $(NSRW_\ee)$ or $(NSOP_\aa)$, we cannot afford to treat separatedly $q$ and $u$: for example in the non-local case, the capillary term is multiplied by a large coefficient $\ee^{-2}$ that would lead to a lifespan of size $\mathcal{O}(\ee^2)$ which is obviously useless for the study of the convergence towards the solution of $(NSK)$. As explained before, as $q$ is not assumed anymore to be small, using the previous a priori estimates on the couple $(q,u)$ leads to an obstruction.
 
The main ingredient in this paper is a new priori estimate in the spirit of \cite{Dtruly} where we recouple $q$ and $u$. Obviously as we need to adapt the idea developped in \cite{Dtruly} where the key is to study equation \eqref{estimub} we cannot hope to use refined estimates obtained from Fourier analysis of the linear system as in \cite{CD, CHVP, Corder}. As we recouple $q$ and $u$, in addition to the term $(1+q)^{-1} \cA u$ we will have to deal with $q \div u$ in the first equation. When $q$ is small this term is harmless and easily absorbed by the left-hand side. On the contrary, in our general case this term prevents any convenient estimates that would be useful for the proof of uniqueness or convergence. It has then to be also included into the linear system we will study. This leads us to the following linear system:
$$
\begin{cases}
\begin{aligned}
&\d_t q+ v\cdot\nabla q+ c\cdot\div u=0,\\
&\d_t u+ v\cdot\nabla u -b\cdot\cA u -\kappa \nabla D[q]=0,\\
\end{aligned}
\end{cases}
$$
where $b,c$ are positive real valued functions, bounded away from zero.

\subsection{Existence}

We refer to the appendix for definitions and properties of the classical and hybrid Besov spaces.

\begin{defi}
\sl{The space $E^s(t)$ is the set of functions $(q,u)$ in
$$
\left(\cC_b([0,t], \dot{B}_{2,1}^s)\cap L_t^1\dot{B}_{2,1}^{s+2}\right) \times
\left(\cC_b([0,t], \dot{B}_{2,1}^{s-1})\cap L_t^1\dot{B}_{2,1}^{s+1}\right)^d
$$
endowed with the norm
\begin{equation}
\|(q,u)\|_{E^s(t)} \overset{def}{=}  \|u\|_{\tilde{L}_t^{\infty} \dot{B}_{2,1}^{s-1}}  +\|q\|_{\tilde{L}_t^{\infty} \dot{B}_{2,1}^{s}} +\|u\|_{L_t^1 \dot{B}_{2,1}^{s+1}} +\|q\|_{L_t^1 \dot{B}_{2,1}^{s+2}}.
\label{normeK}
\end{equation}}
\end{defi}

\begin{defi}
\sl{The space $E_\beta^s(t)$ (for $\beta>0$ expected to be large) is the set of functions $(q,u)$ in
$$
\left(\cC_b([0,t], \dot{B}_{2,1}^s)\cap L_t^1\dot{B}_\beta^{s+2,s}\right) \times
\left(\cC_b([0,t], \dot{B}_{2,1}^{s-1})\cap L_t^1\dot{B}_{2,1}^{s+1}\right)^d
$$
endowed with the norm
\begin{equation}
\|(q,u)\|_{E_\beta^s(t)} \overset{def}{=}  \|u\|_{\tilde{L}_t^{\infty} \dot{B}_{2,1}^{s-1}}  +\|q\|_{\tilde{L}_t^{\infty} \dot{B}_{2,1}^{s}} +\|u\|_{L_t^1 \dot{B}_{2,1}^{s+1}} +\|q\|_{L_t^1 \dot{B}_{\beta}^{s+2,s}}.
\label{normeE}
\end{equation}}
\end{defi}

\begin{rem}
\sl{We observe that the parabolic regularization on $q$ occurs for all frequencies in $E^s(t)$ and only for low frequencies in $E_\beta^s(t)$. Moreover the threshold between the regularized low frequencies and the damped high frequencies goes to infinity as $\beta$ goes to infinity. We refer to \cite{CH, CHVP, Corder} for more details about this threshold and the close relation with the capillary term.}
\end{rem}

\begin{thm}
\sl{Let $\ee>0$, $q_0\in \dot{B}_{2,1}^{\fd}$, $u_0 \in \dot{B}_{2,1}^{\fd-1}$ and assume that $0< \cu \leq 1+q_0 \leq \co$ and $\min(\mu,2\mu+\lambda)>0$. There exist a positive constant $C$ and a time $T>0$, only depending on the physical parameters $d$, $\mu$, $\lambda$, $\kappa$ and the initial data $(q_0,u_0)$, such that system $(NSK)$ has a unique solution $(\rho, u)$ with $(q, u)\in E^{\fd}(T)$, and system $(NSRW_\ee)$ has a unique solution $(\rho_\ee, \ue)$ with $(\qe, \ue)\in E_{1/\ee}^{\fd}(T)$. Moreover
$$
\|(q,u)\|_{E^\fd(T)} +\|(\qe,\ue)\|_{E_{1/\ee}^\fd(T)}  \leq C (\|q_0\|_{\dot{B}_{2,1}^{\fd}} +\|u_0\|_{\dot{B}_{2,1}^{\fd-1}}).
$$
}
\label{thexistR}
\end{thm}

\subsection{Convergence}

As in \cite{CH}, we prove that the solution of $(RW_\ee)$ goes to the solution of $(K)$ when $\ee$ goes to zero.

\begin{thm}
\sl{Assume that $\min(\mu,2\mu+\lambda)>0$, $q_0\in \dot{B}_{2,1}^{\fd}$, $u_0 \in \dot{B}_{2,1}^{\fd-1}$. Then for $T$ given by the previous result,
$$
\|(\qe-q, \ue-u)\|_{E_{1/\ee}^\fd(T)} \underset{\ee \rightarrow 0}{\longrightarrow} 0.
$$
Moreover there exists a constant $C=C(\eta, \kappa, q_0, u_0, T)>0$ such that for all $h\in ]0, 1[$ (if $d=2$) or $h\in ]0,1]$ (if $d\geq 3$),
$$
 \|(\qe-q, \ue-u)\|_{E_{1/\ee}^{\fd-h}(T)}\leq C \ee^h,
$$}
\label{thcvR}
\end{thm}

\subsection{Order parameter model}

The very same results are true for the order parameter model:

\begin{defi}
\sl{The space $F_\beta^s(t)$ (for $\beta>0$) is the set of functions $(q,c,u)$ in
$$
\left(\cC_b([0,t], \dot{B}_{2,1}^s)\cap L_t^1\dot{B}_\beta^{s+2,s}\right)^2 \times
\left(\cC_b([0,t], \dot{B}_{2,1}^{s-1})\cap L_t^1\dot{B}_{2,1}^{s+1}\right)^d
$$
endowed with the norm
\begin{multline}
\|(q,c,u)\|_{F_\beta^s(t)} \overset{def}{=}  \|u\|_{\tilde{L}_t^{\infty} \dot{B}_{2,1}^{s-1}}  +\|q\|_{\tilde{L}_t^{\infty} \dot{B}_{2,1}^{s}} +\|c\|_{\tilde{L}_t^{\infty} \dot{B}_{2,1}^{s}}\\
 +\|u\|_{L_t^1 \dot{B}_{2,1}^{s+1}} +\|q\|_{L_t^1 \dot{B}_{\beta}^{s+2,s}} +\|c\|_{L_t^1 \dot{B}_{\beta}^{s+2,s}}.
\label{normeF}
\end{multline}}
\end{defi}

\begin{thm}
\sl{Let $\aa>0$, $q_0\in \dot{B}_{2,1}^{\fd}$, $u_0 \in \dot{B}_{2,1}^{\fd-1}$ and assume $0< \cu \leq 1+q_0 \leq \co$ and $\min(\mu,2\mu+\lambda)>0$. Let $c_0$ be defined by $-\Delta c_0+ \aa^2 c_0=\aa^2 \rho_0$, that is $c_0=\psia *\rho_0$. There exist a positive constant $C$  and a time $T>0$ only depending on the physical parameters and $(q_0,u_0)$ such that system $(NSK)$ has a unique solution $(\rho, u)$ with $(q, u)\in E^{\fd}(T)$, and system $(NSOP_\aa)$ has a unique solution $(\rho_\aa, \ca, \ua)$ with $(\qa, \ca,\ua)\in F_{\aa}^{\fd}(T)$ and $\ca= \psia*\qa$. Moreover:
$$
\|(q,u)\|_{E^\fd(T)} +\|(\qa, \ca, \ua)\|_{F_{\aa}^\fd(T)}  \leq C (\|q_0\|_{\dot{B}_{2,1}^{\fd}} +\|u_0\|_{\dot{B}_{2,1}^{\fd-1}}),
$$
and
$$
\|\ca-\rho_\aa\|_{\tilde{L}^\infty_T\dot{B}_{2,1}^{\fd}} \underset{\aa\rightarrow \infty}{\longrightarrow} 0,\mbox{ and } \|\ca-\rho_\aa\|_{L^1_T\dot{B}_{2,1}^{\fd}} \leq C \aa^{-2}.
$$
}
\label{thexistOP}
\end{thm}

\begin{thm}
\sl{Assume that $\min(\mu,2\mu+\lambda)>0$, $q_0\in \dot{B}_{2,1}^{\fd}$, $u_0 \in \dot{B}_{2,1}^{\fd-1}$. Then for $T$ given by the previous result,
$$
\|(\qa-q, c_\aa-\rho, \ua-u)\|_{F_{\aa}^{\fd}(T)} \underset{\aa \rightarrow \infty}{\longrightarrow} 0.
$$
Moreover there exists a constant $C=C(\eta, \kappa, q_0, T)>0$ such that for all $h\in ]0, 1[$ (if $d=2$) or $h\in ]0,1]$ (if $d\geq 3$), and for all $t\in[0,T]$,
$$
 \|(\qa-q, c_\aa-\rho, \ua-u)\|_{F_{\aa}^{\fd-h}(T)}\leq C \aa^{-h},
$$}
\label{thcvOP}
\end{thm}

We refer to the end of the article for the variable coefficients case, and the case of Besov spaces $\dot{B}_{p,1}^s$ with $p\neq 2$.

\section{Proof of Theorem \ref{thexistR}}

We will prove theorems \ref{thexistR} and \ref{thcvR}. As we only use energy methods, the proofs are strictly the same for the order parameter model, because in the $L^2$-setting, the interaction potentials $\psia$ and $\phie$ play exactly the same role. Let us emphasize that this was not the case in \cite{CHVP, Corder} as it involved Fourier computations and finite differences representations of Besov norms.

\subsection{Linear estimates with variable coefficients}

As announced, these results rely on a priori estimates for solutions of the following system:
\begin{equation}
\begin{cases}
\begin{aligned}
&\d_t q+ v\cdot\nabla q+ c\cdot\div u= F,\\
&\d_t u+ v\cdot\nabla u -b\cdot\cA u-\kappa \nabla D[q]= G.\\
\end{aligned}
\end{cases}
\label{SystL}
\end{equation}
With
$$\cA u= \mu \Delta u+ (\lambda+\mu)\nabla \div u,$$
and we recall that $D[q]$ is given in \eqref{Capill}.
\begin{thm}
\sl{Let $s\in]-\fd,\fd+1]$, $\nu=2\mu+\lambda$ and assume that $\nu_0=\min(\mu,\nu)>0$, $q_0\in \dot{B}_{2,1}^{s}$, $u_0 \in \dot{B}_{2,1}^{s-1}$. Let $T>0$, and assume that $(q,u)$ solves system \eqref{SystL} on $[0,T]$, with $v\in L_T^\infty \dot{B}_{2,1}^{\fd-1} \cap L_T^1 \dot{B}_{2,1}^{\fd+1}$ and $b,c:[0,T]\times \R^d \rightarrow \R_+$ such that:
\begin{itemize}
\item $b-1,c-1 \in \mathcal{C}([0,T], \dot{B}_{2,1}^\fd),$
\item $\d_t b, \d_t c \in L_T^2\dot{B}_{2,1}^{\fd-1},$
\item for all $t\leq T, x\in\R^d$, we have $0<c_* \leq c(t,x) \leq c^*$ and $0<b_* \leq b(t,x) \leq b^*.$
\end{itemize}
Then $(q,u)\in E^s(T)$ (respectively $(\qe, \ue) \in E_{1/\ee}^s(T)$, $(\qa, \ua) \in E_{\aa}^s(T)$). Moreover if we define
$$
g^s(q,u)(t)= \Sum_{j\in \Z} 2^{j(s-1)} \sup_{t' \in [0,t]} \left(\|u_j(t')\|_{L^2} + h_j(t')\right),
$$
with
\begin{multline}
h_j(t')^2 =(u_j(t')|c_m u_j(t'))_{L^2} +\kappa (q_j(t')|D[q_j(t')])_{L^2}\\
+\eta \left(2 (u_j(t')|\nabla q_j(t'))_{L^2}+ \nu (\nabla q_j(t')|\frac{b_m}{c_m} \nabla q_j(t'))_{L^2}  \right),
\end{multline}
where $u_j=\ddj u$ (the same for $q$) and $b_m, c_m$ are smooth functions defined by:
\begin{equation}
b_m=1+\dot{S}_m (b-1) \quad \mbox{ and } \quad c_m=1+\dot{S}_m (c-1).
\label{defbmcm}
\end{equation}
Then there exist $m_0 \in \Z$, two constants $\gamma_*>0$ and $\Fs\geq 1$ such that if $\eta>0$ is fixed small enough (all of them only depending on the bounds $b_*,c_*,b^*,c^*$ and the viscous and capillary coefficients) then for all $t\leq T$ and $m\geq m_0$,
\begin{equation}
{\Fs}^{-1} g^s(t) \leq \|u\|_{\tilde{L}_t^\infty \dot{B}_{2,1}^{s-1}} +\|q\|_{\tilde{L}_t^\infty \dot{B}_{2,1}^s} \leq \Fs g^s(t),
\label{equivgs}
\end{equation}
and
\begin{multline}
g^s(q,u)(t)+\frac{\nu_0 b_*}{4} \|u\|_{L_t^1 \dot{B}_{2,1}^{s+1}} +\gamma_* \|D[q]\|_{L_t^1 \dot{B}_{2,1}^{s}} \leq g^s(q_0, u_0)(0) +\Fs \int_0^t \left(\|F\|_{\dot{B}_{2,1}^s} +\|G\|_{\dot{B}_{2,1}^{s-1}}\right) d\tau\\
+ \Fs \int_0^t g^s(q,u)(\tau) \Bigg[ 2^m \Big(\|\d_t b\|_{\dot{B}_{2,1}^{\fd-1}} +\|\d_t c\|_{\dot{B}_{2,1}^{\fd-1}}\Big)\\
+(1+\|b-1\|_{\dot{B}_{2,1}^\fd} +\|c-1\|_{\dot{B}_{2,1}^\fd})^2 \Big(\|\nabla v\|_{\dot{B}_{2,1}^\fd} +2^m \|v\|_{\dot{B}_{2,1}^\fd} +2^{2m} +\|v\|_{\dot{B}_{2,1}^\fd}^2\Big)\Bigg] d\tau\\
+\Fs \int_0^t \left(\|(I_d-\dot{S}_m)(b-1)\|_{\dot{B}_{2,1}^\fd} +\|(I_d-\dot{S}_m)(c-1)\|_{\dot{B}_{2,1}^\fd}\right) \|u\|_{\dot{B}_{2,1}^{s+1}} d\tau.
\label{estimgs}
\end{multline}
}
\label{thestimbc}
\end{thm}
\begin{rem}
\sl{The condition $s\in]-\fd,\fd+1]$ is required by paraproduct and remainder laws, we refer to \eqref{Bony} for details.}
\end{rem}

\begin{rem}
\sl{Let us emphasize that we have:
$$
\|D[q]\|_{L_T^1 \dot{B}_{2,1}^{s}} \sim
\begin{cases}
\|q\|_{L_T^1 \dot{B}_{2,1}^{s+2}} \mbox{ in the local case }(NSK),\\
\|q\|_{L_T^1 \dot{B}_{1/\ee}^{s+2,s}} \mbox{ in the non-local case }(NSRW_\ee),\\
\|q\|_{L_T^1 \dot{B}_{\aa}^{s+2,s}} \mbox{ in the non-local case }(NSOP_\aa).
\end{cases}
$$
}
\end{rem}

\subsection{Proof of Theorem \ref{thestimbc}}

We will prove the result in the first non-local case, that is for $D[q]=\displaystyle{\frac{\phie *q-q}{\ee^2}}$. As said before, for the order parameter model everything works the same, and for the local case the same argument is valid but many steps are much easier thanks to the fact that $q$ is more regular. We will highlight in the following proof what can be simplified for the local case $(NSK)$.

From the definition of $b_m$ we have $b_m-1=\dot{S}_m (c-1)$, so that $b_m-1$ is smooth and its $\dot{B}_{2,1}^\fd$-norm is bounded by the one of $b-1$. Moreover $b-b_m=(I_d-\dot{S}_m)(b-1)$, and we can fix $m_0$ large enough so that for all $m\geq m_0$:
$$
\|b-b_m\|_{\dot{B}_{2,1}^\fd}\leq \frac{b_*}{2} \mbox{ and }\|c-c_m\|_{\dot{B}_{2,1}^\fd}\leq \frac{c_*}{2}.
$$
In this case, thanks to the injection $\dot{B}_{2,1}^\fd \hookrightarrow L^\infty$, we immediately have for all $t\leq T, x\in\R^d$, that
\begin{equation}
0<\frac{c_*}{2} \leq c_m(t,x) \leq c^* +\frac{c_*}{2}\quad \mbox{ and } \quad 0<\frac{b_*}{2}  \leq b_m(t,x) \leq b^*+\frac{b_*}{2}.
\label{bornesbc}
\end{equation}
Next, in the spirit of \cite{Dtruly}, let us first rewrite system \eqref{SystL} as follows:
$$
\begin{cases}
\begin{aligned}
&\d_t q+ v\cdot\nabla q+ c_m\div u= F+F_m,\\
&\d_t u+ v\cdot\nabla u -b_m\cA u-\kappa \frac{\phie*\nabla q-\nabla q}{\ee^2}= G+G_m,\\
\end{aligned}
\end{cases}
$$
where
$$
\begin{cases}
F_m=(c_m-c)\cdot \div u =-\left(I_d-\dot{S}_m\right)(c-1) \cdot \div u,\\
G_m=(b-b_m) \cdot \cA u =\left(I_d-\dot{S}_m\right)(b-1) \cdot \cA u.
\end{cases}
$$
Then applying operator $\ddj$ to the system, and using the notation $f_j=\ddj f$, we obtain (following the lines of \cite{Dtruly}):
\begin{equation}
\begin{cases}
\begin{aligned}
&\d_t q_j+ v\cdot\nabla q_j+\div(c_m u_j) =f_j,\\
&\d_t u_j+ v\cdot\nabla u_j -\mu \div(b_m.\nabla u_j)-(\lambda+\mu) \nabla (b_m \div u_j)-\kappa \frac{\phie *\nabla q_j-\nabla q_j}{\ee^2} =g_j,
\end{aligned}
\end{cases}
\label{systloc}
\end{equation}
where
\begin{equation}
\begin{cases}
f_j= F_j+F_{m,j}+R_j+\tilde{R}_j \quad \quad \mbox{and} \quad g_j=G_j+G_{m,j}+S_j+\tilde{S}_j,\\
R_j=v\cdot q_j-\ddj (v\cdot q) \quad \quad \mbox{and} \quad S_j=v\cdot u_j-\ddj (v\cdot u),\\
\tilde{R}_j = \div (c_m u_j)-\ddj(c_m \div u) =\div \Big((c_m-1) u_j\Big)-\ddj\Big((c_m-1) \div u\Big).\\
\tilde{S}_j = \mu\left(\ddj (b_m \Delta u)-\div(b_m \nabla \ddj u) \right)+ (\lambda+\mu) \left( \ddj (b_m \nabla \div u)- \nabla (b_m \div \ddj u)\right).\\
\hspace{0.5cm}
=\mu\left(\ddj ((b_m-1) \Delta u)-\div((b_m-1) \nabla \ddj u) \right)\\
\hspace{4cm}
+ (\lambda+\mu) \left( \ddj ((b_m-1) \nabla \div u)- \nabla ((b_m-1) \div \ddj u)\right).
\end{cases}
\label{deffext}
\end{equation}
Let us begin by stating estimates on these external terms (we refer to lemma \ref{estimtruly} in the appendix and \cite{Dtruly, Dbook} for details and proofs):
\begin{prop}
\sl{Under the previous assumptions, there exist a positive constant $C$ and a nonnegative summable sequence $(c_j)_{j\in \Z} =\big(c_j(t)\big)_{j\in \Z}$ whose summation is $1$ such that if we denote by $\overline{\nu}=\mu +|\lambda +\mu|$, then for all $j\in \Z$ we have:
\begin{equation}
\begin{cases}
\|R_j\|_{L^2} \leq C c_j 2^{-js} \|\nabla v\|_{\dot{B}_{2,1}^\fd} \|q\|_{\dot{B}_{2,1}^s},\\
\|S_j\|_{L^2} \leq C c_j 2^{-j(s-1)} \|\nabla v\|_{\dot{B}_{2,1}^\fd} \|u\|_{\dot{B}_{2,1}^{s-1}},\\
\|\tilde{R}_j\|_{L^2} \leq C c_j 2^{-js} 2^m \|c-1\|_{\dot{B}_{2,1}^\fd} \|u\|_{\dot{B}_{2,1}^s},\\
\|\tilde{S}_j\|_{L^2} \leq C \overline{\nu} c_j 2^{-j(s-1)} 2^m \|b-1\|_{\dot{B}_{2,1}^\fd} \|u\|_{\dot{B}_{2,1}^s},\\
\|F_m\|_{\dot{B}_{2,1}^s} \leq C \|\left(I_d-\dot{S}_m \right)(c-1)\|_{\dot{B}_{2,1}^\fd} \|u\|_{\dot{B}_{2,1}^{s+1}},\\
\|G_m\|_{\dot{B}_{2,1}^{s-1}} \leq C \|\left(I_d-\dot{S}_m \right)(b-1)\|_{\dot{B}_{2,1}^\fd} \|u\|_{\dot{B}_{2,1}^{s+1}}.
\end{cases}
\label{estimfext}
\end{equation}
}
\end{prop}
In this proof of theorem \ref{thestimbc}, the ideas are classical (we refer to \cite{Dinv, Dtruly, Dheatlocal, Dbook, Has2, Has5, CH}) and consist of combining innerproducts in $L^2$ of the equations in order to cancel terms that we are not able to estimate (too much derivatives or large coefficients).
\begin{rem}
\sl{Due to the initial regularity, the most natural way is to study the evolution of $\|u_j\|_{L^2}^2$ and $\|\nabla q_j\|_{L^2}^2$, for this we consider the inner product of the gradient of the first equation by $\nabla q_j$, and the second by $u_j$. This computation involves the following problematic term $\kappa \ee^{-2} (\phie*\nabla q_j-\nabla q_j|u_j)_{L^2}$ that unfortunately, at this level of the study, we are not able to estimate uniformly with respect to $\ee$. We need a way to cancel it and the easiest way to do this is, as in \cite{CH}, to consider the innerproduct of the density equation by $\ee^{-2}(q_j-\phie *q_j)$ and the velocity equation by $c_m u_j$. Then the problematic term will be neutralized if we sum the results.}
\label{rqintegraleepsilon}
\end{rem}
Keeping in mind the fact that $\frac{d}{dt} (u_j| c_m u_j)_{L^2}= 2(\d_t u_j| c_m u_j)_{L^2}+( \d_t c_m.u_j|u_j)_{L^2}$, we begin by taking the inner product of the velocity equation with $c_m u_j$:
\begin{multline}
(\d_t u_j| c_m u_j)_{L^2}+ (v\cdot\nabla u_j| c_m u_j)_{L^2} +\mu (b_m.\nabla u_j| \nabla(c_m u_j))_{L^2}+(\lambda+\mu) (b_m \div u_j| \div (c_m u_j))_{L^2}\\
-\kappa (\frac{\phie *\nabla q_j-\nabla q_j}{\ee^2}| c_m u_j)_{L^2} =(g_j| c_m u_j)_{L^2}.
\end{multline}
We have:
$$
\begin{cases}
\displaystyle{(b_m.\nabla u_j| \nabla(c_m u_j))_{L^2} =\int_{\R^d} b_m c_m.|\nabla u_j|^2 dx +(b_m.\nabla u_j| u_j.\nabla c_m)_{L^2}},\\
\displaystyle{(b_m \div u_j| \div (c_m u_j))_{L^2} =\int_{\R^d} b_m c_m.|\div u_j|^2 dx +(b_m \div u_j|u_j.\nabla c_m)_{L^2}},
\end{cases}
$$
with
$$
\Big| \mu (b_m.\nabla u_j| u_j.\nabla c_m)_{L^2} +(\lambda+\mu) (b_m \div u_j|u_j.\nabla c_m)_{L^2}\Big| \leq \overline{\nu} \|b_m\|_{L^\infty} \|\nabla c_m\|_{L^\infty} 2^j \|u_j\|_{L^2}^2.
$$
We estimate the following term like in \cite{Dinv, CH} using integrations by parts:
\begin{multline}
\Big|(v\cdot\nabla u_j| c_m u_j)_{L^2}\Big| \leq \frac{1}{2} \|\div (c_m v)\|_{L^\infty} \|u_j\|_{L^2}^2\\
\leq \Big(\|c_m\|_{L^\infty}.\|\nabla v\|_{L^\infty}+\|v\|_{L^\infty}.\|\nabla c_m\|_{L^\infty} \Big) \|u_j\|_{L^2}^2.
\end{multline}
Thanks to the frequency truncation in the definition of $c_m$ and the fact that
$$
\d_t c_m =\d_t (c_m-1)= \d_t \left(\dot{S}_m(c-1)\right) =\dot{S}_m \left(\d_t(c-1)\right) =\dot{S}_m \left(\d_t c\right),$$
the additional term is estimated by:
$$
( \d_t c_m.u_j|u_j)_{L^2} \leq \|\d_t c_m\|_{L^\infty} \|u_j\|_{L^2}^2 \leq \|\d_t c_m\|_{\dot{B}_{2,1}^\fd} \|u_j\|_{L^2}^2 \leq 2^m \|\d_t c\|_{\dot{B}_{2,1}^{\fd-1}} \|u_j\|_{L^2}^2.
$$
Gathering these estimates we obtain:
\begin{multline}
\frac{1}{2} \frac{d}{dt}(u_j| c_m u_j)_{L^2} +\mu \int_{\R^d} b_m c_m.|\nabla u_j|^2 dx +(\lambda+\mu) \int_{\R^d} b_m c_m.|\div u_j|^2 dx\\
-\kappa (\frac{\phie *\nabla q_j-\nabla q_j}{\ee^2}| c_m u_j)_{L^2} \leq \|g_j\|_{L^2} \|c_m\|_{L^\infty} \|u_j\|_{L^2} +2^m \|\d_t c\|_{\dot{B}_{2,1}^{\fd-1}} \|u_j\|_{L^2}^2\\
+\overline{\nu} \|b_m\|_{L^\infty} \|\nabla c_m\|_{L^\infty} 2^j \|u_j\|_{L^2}^2 +\Big(\|c_m\|_{L^\infty}.\|\nabla v\|_{L^\infty}+\|v\|_{L^\infty}.\|\nabla c_m\|_{L^\infty} \Big) \|u_j\|_{L^2}^2.
\end{multline}
Thanks to the bounds on $b_m$ and $c_m$ (see \eqref{bornesbc}), we prove similarly to \cite{Dinv, CH} that (wether $\lambda+\mu$ is negative or not) that:
$$
\mu \int_{\R^d} b_m c_m.|\nabla u_j|^2 dx +(\lambda+\mu) \int_{\R^d} b_m c_m.|\div u_j|^2 dx \geq \nu_0 \frac{b_* c_*}{4} 2^{2j} \|u_j\|_{L^2}^2,
$$
where we recall that
$$
\begin{cases}
\nu_0=\min(\mu, 2\nu+\lambda),\\
\overline{\nu}=\mu +|\mu+\lambda|.
\end{cases}
$$
Moreover,
$$
\overline{\nu} \|b_m\|_{L^\infty} \|\nabla c_m\|_{L^\infty} 2^j \|u_j\|_{L^2}^2 \leq \nu_0 \frac{b_* c_*}{8} 2^{2j} \|u_j\|_{L^2}^2 +\frac{2 \overline{\nu}^2}{\nu_0 b_* c_*} \|b_m\|_{L^\infty}^2 \|\nabla c_m\|_{L^\infty}^2 \|u_j\|_{L^2}^2.
$$
Using the bound from \eqref{bornesbc}, and the fact that
$$
\|\nabla c_m\|_{L^\infty} \leq \|\nabla (c_m-1)\|_{\dot{B}_{2,1}^\fd} \leq \|\nabla \dot{S}_m(c-1)\|_{\dot{B}_{2,1}^\fd} \leq 2^m \|c-1\|_{\dot{B}_{2,1}^\fd},
$$
we end up with:
\begin{multline}
\frac{1}{2} \frac{d}{dt}(u_j| c_m u_j)_{L^2} +\nu_0 \frac{b_* c_*}{8} 2^{2j} \|u_j\|_{L^2}^2 -\kappa (\frac{\phie *\nabla q_j-\nabla q_j}{\ee^2}| c_m u_j)_{L^2}\\
\leq \Cs \|g_j\|_{L^2} \|u_j\|_{L^2} + C \|u_j\|_{L^2}^2 \Bigg[ 2^m \|\d_t c\|_{\dot{B}_{2,1}^{\fd-1}} +\Cs \Big(\|\nabla v\|_{L^\infty}+ 2^m \|c-1\|_{\dot{B}_{2,1}^\fd} \|v\|_{L^\infty} \Big)\\
+\Cs \frac{\overline{\nu}^2}{\nu_0} 2^{2m} \|c-1\|_{\dot{B}_{2,1}^\fd}^2 \Bigg].
\label{estimuj}
\end{multline}
where $\Cs\geq 1$ denotes a constant only depending on the bounds of $b$ and $c$
\begin{rem}
\sl{In the following we adopt the convention that even if varying from line to line we will always denote by $\Cs$ such a constant, and by $\Ds, \Fs$ a constant that in addition depends on the physical parameters ($\lambda$, $\mu$, $\kappa$).}
\end{rem}
Let us now turn to the density fluctuation: as explained, in order to neutralize the capillary term, instead of studying $\|\nabla q_j\|^2$, we study the quantity $(q_j| \frac{q_j-\phie *q_j}{\ee^2})_{L^2}$. computing the inner product of the density equation by $\frac{q_j-\phie *q_j}{\ee^2}$, we obtain that:
\begin{multline}
\frac{1}{2} \frac{d}{dt} (q_j| \frac{q_j-\phie *q_j}{\ee^2})_{L^2} +(v\cdot\nabla q_j| \frac{q_j-\phie *q_j}{\ee^2})_{L^2} +(\div(c_m u_j)| \frac{q_j-\phie *q_j}{\ee^2})_{L^2}\\
=(f_j| \frac{q_j-\phie *q_j}{\ee^2})_{L^2}.
\end{multline}
We refer to \cite{CHVP, Corder} and the appendix for the fact that:
$$
\begin{cases}
\displaystyle{(q_j| \frac{q_j-\phie *q_j}{\ee^2})_{L^2} \sim \min (\frac{1}{\ee^2}, 2^{2j}) \|q_j\|_{L^2}^2,}\\
\displaystyle{\|\frac{q_j-\phie *q_j}{\ee^2}\|_{L^2} \sim \min (\frac{1}{\ee^2}, 2^{2j}) \|q_j\|_{L^2}.}
\end{cases}
$$
So that, exactly like in \cite{CH}, we estimate:
\begin{multline}
\Big|(f_j| \frac{q_j-\phie *q_j}{\ee^2})_{L^2}\Big| \leq \|f_j\|_{L^2} \|\frac{q_j-\phie *q_j}{\ee^2}\|_{L^2} \leq \|f_j\|_{L^2} \min (\frac{1}{\ee^2}, 2^{2j}) \|q_j\|_{L^2}\\
\leq 2^j \|f_j\|_{L^2} \sqrt{\min (\frac{1}{\ee^2}, 2^{2j}) \|q_j\|_{L^2}^2} \leq C 2^j \|f_j\|_{L^2} \sqrt{(q_j| \frac{q_j-\phie *q_j}{\ee^2})_{L^2}},
\end{multline}
and similarly,
$$
\Big|(v\cdot\nabla q_j| \frac{q_j-\phie *q_j}{\ee^2})_{L^2}\Big| \leq C\|v\|_{L^\infty} \|\nabla q_j\|_{L^2} 2^j \sqrt{(q_j| \frac{q_j -\phie *q_j}{\ee^2})_{L^2}}.
$$
Collecting these estimates implies that:
\begin{multline}
\frac{1}{2} \frac{d}{dt} (q_j| \frac{q_j-\phie *q_j}{\ee^2})_{L^2} -(c_m u_j| \frac{\nabla q_j-\phie *\nabla q_j}{\ee^2})_{L^2}\\
\leq C 2^j \|f_j\|_{L^2} \sqrt{(q_j| \frac{q_j-\phie *q_j}{\ee^2})_{L^2}} +\|v\|_{L^\infty} \|\nabla q_j\|_{L^2} 2^j \sqrt{(q_j| \frac{q_j -\phie *q_j}{\ee^2})_{L^2}}.
\label{estimqj}
\end{multline}
so that, when we compute \eqref{estimuj}$+\kappa$\eqref{estimqj}, there is a cancellation of the problematic terms $\ee^{-2}(\div(c_m u_j)| q_j-\phie *q_j)_{L^2}$ and we finally obtain:
\begin{multline}
\frac{1}{2} \frac{d}{dt}\left( (u_j| c_m u_j)_{L^2} +\kappa (q_j| \frac{q_j-\phie *q_j}{\ee^2})_{L^2} \right) +\nu_0 \frac{b_* c_*}{8} 2^{2j} \|u_j\|_{L^2}^2 \leq \Cs \|g_j\|_{L^2} \|u_j\|_{L^2} \\
+ C\kappa 2^j \|f_j\|_{L^2} \sqrt{(q_j| \frac{q_j-\phie *q_j}{\ee^2})_{L^2}} +C \kappa \|v\|_{L^\infty} \|\nabla q_j\|_{L^2} 2^j \sqrt{(q_j| \frac{q_j -\phie *q_j}{\ee^2})_{L^2}}\\
+ C \|u_j\|_{L^2}^2 \Bigg[ 2^m \|\d_t c\|_{\dot{B}_{2,1}^{\fd-1}} +\Cs \Big(\|\nabla v\|_{L^\infty}+ 2^m \|c-1\|_{\dot{B}_{2,1}^\fd} \|v\|_{L^\infty} \Big) +\Cs \frac{\overline{\nu}^2}{\nu_0} 2^{2m} \|c-1\|_{\dot{B}_{2,1}^\fd}^2 \Bigg].
\label{estimujqj1}
\end{multline}
As in \cite{CH}, the only way to obtain a regularization for the density fluctuation is to consider the inner product of the velocity equation by $\nabla q_j$ which provides the nonnegative term $\kappa \ee^{-2}(\nabla q_j| \nabla q_j-\phie * \nabla q_j)_{L^2}$, that is we will study the variation in time of $(u_j|\nabla q_j)_{L^2}$. As
$$
\frac{d}{dt} (u_j|\nabla q_j)_{L^2}=  (\d_t u_j|\nabla q_j)_{L^2} + (u_j|\d_t \nabla q_j)_{L^2},
$$
like in \cite{CH} we are lead to sum the following estimates:
\begin{multline}
(\d_t u_j| \nabla q_j)_{L^2}+ (v\cdot\nabla u_j| \nabla q_j)_{L^2} -\mu (\div(b_m.\nabla u_j)| \nabla q_j)_{L^2}-(\lambda+\mu) (\nabla(b_m \div u_j)| \nabla q_j)_{L^2}\\
-\kappa (\nabla q_j|\frac{\phie *\nabla q_j-\nabla q_j}{\ee^2})_{L^2} =(g_j| \nabla q_j)_{L^2} \leq \|g_j\|_{L^2} \|\nabla q_j\|_{L^2},
\end{multline}
and (taking the innerproduct of the equation on $\nabla q_j$ with $u_j$)
\begin{multline}
(\d_t \nabla q_j|u_j)_{L^2} +(\nabla(v\cdot \nabla q_j)|u_j)_{L^2} -(\div (c_m u_j)|\div u_j)_{L^2}=(\nabla f_j|u_j)_{L^2}=-(f_j|\div u_j)_{L^2}\\ \leq \|f_j\|_{L^2} 2^j \|u_j\|_{L^2}.
\end{multline}
A simple computation shows that (we refer to \cite{Dinv, Dbook, CH} for details):
$$
\Big|(v\cdot\nabla u_j| \nabla q_j)_{L^2} +(\nabla(v\cdot \nabla q_j)|u_j)_{L^2} \Big| \leq C \|\nabla v\|_{L^\infty} \|\nabla q_j\|_{L^2} \|u_j\|_{L^2}.
$$
Moreover, a rough estimate on the term
$$
(\div (c_m u_j)|\div u_j)_{L^2} =(c_m \div u_j +u_j \cdot \nabla c_m|\div u_j)_{L^2},
$$
gives, for the same reasons as previously
\begin{multline}
\Big| (\div (c_m u_j)|\div u_j)_{L^2} \Big| \leq \|c_m\|_{L^\infty} \|\nabla u_j\|_{L^2}^2 +\|\nabla c_m\|_{L^\infty} 2^j \|u_j\|_{L^2}^2\\
\leq \|c_m\|_{L^\infty} 2^{2j}\| u_j\|_{L^2}^2 +2^m\|c-1\|_{\dot{B}_{2,1}^\fd} 2^j \|u_j\|_{L^2}^2\\
\leq \Cs 2^{2j}\|u_j\|_{L^2}^2 +2^{2m}\|c-1\|_{\dot{B}_{2,1}^\fd}^2 \|u_j\|_{L^2}^2,
\end{multline}
Where, as usual, $\Cs$ is a constant only depending on the bounds of $b$ and $c$. Finally:
\begin{multline}
\frac{d}{dt} (u_j|\nabla q_j)_{L^2} -\mu (\div(b_m.\nabla u_j)| \nabla q_j)_{L^2}-(\lambda+\mu) (\nabla(b_m \div u_j)| \nabla q_j)_{L^2}\\
+\kappa (\nabla q_j|\frac{\nabla q_j-\phie *\nabla q_j}{\ee^2})_{L^2} \leq \|g_j\|_{L^2} \|\nabla q_j\|_{L^2} +2^j \|f_j\|_{L^2} \|u_j\|_{L^2}\\
+C \|\nabla v\|_{L^\infty} \|\nabla q_j\|_{L^2} \|u_j\|_{L^2} +\Cs 2^{2j}\|u_j\|_{L^2}^2 +2^{2m}\|c-1\|_{\dot{B}_{2,1}^\fd}^2 \|u_j\|_{L^2}^2.
\label{estimujgradqj}
\end{multline}
This estimate provides regularization for $(q_j| \frac{q_j-\phie *q_j}{\ee^2})_{L^2}$ (this is the best we can hope for), but involves terms that cannot be absorbed or neutralized through the Gronwall lemma because they introduce too many derivatives, namely:
\begin{equation}
-\mu (\div(b_m.\nabla u_j)| \nabla q_j)_{L^2}-(\lambda+\mu) (\nabla(b_m \div u_j)| \nabla q_j)_{L^2}.
\label{difftermes}
\end{equation}
As in \cite{CH} we need to compensate them thanks to the density equation: more precisely in \cite{CH} (constant coefficients, no $b_m$) we could entirely neutralize $-\mu (\Delta u_j)| \nabla q_j)_{L^2}-(\lambda+\mu) (\nabla\div u_j| \nabla q_j)_{L^2}$ with $(\nabla\div u_j| \nabla q_j)_{L^2}$ from the estimate on $\|\nabla q_j\|_{L^2}^2$ thanks to integrations by parts. In our case, due to the variable coefficients $b_m$ and $c_m$, it will not be possible to entirely cancel those terms, but we will be able to substract to \eqref{difftermes} its most dangerous parts, that is where all the derivatives pound on $u_j$ or $q_j$, the rest being absorbable because at least one derivative pounds on $b_m$ or $c_m$ (and then producing a harmless $2^m$ thanks to the Bernstein lemma). For this, we study the variations of $(\nabla q_j|\frac{b_m}{c_m}\nabla q_j)_{L^2}$. As we have
\begin{equation}
\frac{d}{dt} (\nabla q_j|\frac{b_m}{c_m}\nabla q_j)_{L^2} = 2(\d_t \nabla q_j|\frac{b_m}{c_m}\nabla q_j)_{L^2} +(\nabla q_j|\d_t(\frac{b_m}{c_m})\nabla q_j)_{L^2},
\label{estimddt}
\end{equation}
we begin by estimating the derivative of $\d_t(\frac{b_m}{c_m})$: there exists a constant only depending on the bounds of $b$ and $c$ once again denoted by $\Cs$ so that:
\begin{multline}
\|\d_t(\frac{b_m}{c_m})\|_{L^\infty} = \|\frac{\d_t b_m}{c_m}-\frac{b_m}{c_m^2} \d_t c_m\|_{L^\infty} \leq \frac{\|\d_t b_m\|_{\dot{B}_{2,1}^\fd}}{\frac{c_*}{2}}+ \frac{\|b_m\|_{L^\infty}}{(\frac{c_*}{2})^2}\|\d_t c_m\|_{\dot{B}_{2,1}^\fd}\\
\leq \Cs (\|\d_t b_m\|_{\dot{B}_{2,1}^\fd}+ \|\d_t c_m\|_{\dot{B}_{2,1}^\fd}) \leq 2^m \Cs (\|\d_t b\|_{\dot{B}_{2,1}^{\fd-1}}+ \|\d_t c\|_{\dot{B}_{2,1}^{\fd-1}}).
\label{estimdtbmcm}
\end{multline}
Next, taking the inner product of the equation on $\nabla q_j$ by $\frac{b_m}{c_m}\nabla q_j$,
$$
(\d_t \nabla q_j|\frac{b_m}{c_m}\nabla q_j)_{L^2}+ (\nabla(v\cdot\nabla q_j)|\frac{b_m}{c_m}\nabla q_j)_{L^2} +(\nabla \div(c_m u_j)|\frac{b_m}{c_m}\nabla q_j)_{L^2} =(\nabla f_j|\frac{b_m}{c_m}\nabla q_j)_{L^2}.
$$
The last term of the left-hand side will help us to balance the problematic terms from before, and we begin by estimating the other terms. Without any surprise, as in \cite{Dinv, Dbook, CH}, thanks to integrations by parts, we have:
\begin{multline}
(\nabla(v\cdot\nabla q_j)|\frac{b_m}{c_m}\nabla q_j)_{L^2} = \int_{\R^d} \Sum_{k,l=1}^d \left( \d_k v^l.\d_l q_j.\frac{b_m}{c_m} \d_k q_j +v^l \frac{b_m}{c_m} \frac{1}{2} \d_l(\d_k q_j)^2 \right) dx\\
= \int_{\R^d} \left(\Sum_{k,l=1}^d \d_k v^l.\d_l q_j.\frac{b_m}{c_m} \d_k q_j  -\frac{1}{2} \div(\frac{b_m}{c_m}v)  |\nabla q_j|^2 \right)dx,
\end{multline}
so that
\begin{multline}
\Big|(\nabla(v\cdot\nabla q_j)|\frac{b_m}{c_m}\nabla q_j)_{L^2}\Big| \leq C \left(\|\nabla v\|_{L^\infty} \|\frac{b_m}{c_m}\|_{L^\infty} +\|\div(\frac{b_m}{c_m}v)\|_{L^\infty}\right) \|\nabla q_j\|_{L^2}^2\\
\leq C \left(2 \|\nabla v\|_{L^\infty} \|\frac{b_m}{c_m}\|_{L^\infty} +\|v\|_{L^\infty} \|\nabla(\frac{b_m}{c_m})\|_{L^\infty}\right) \|\nabla q_j\|_{L^2}^2\\
\end{multline}
Similarly to \eqref{estimdtbmcm}, there exists a constant $\Cs$ only depending on the bounds of $b$ and $c$ so that:
\begin{multline}
\|\nabla(\frac{b_m}{c_m})\|_{L^\infty} = \|\frac{\nabla b_m}{c_m}-\frac{b_m}{c_m^2} \nabla c_m\|_{L^\infty} \leq \frac{\|\nabla b_m\|_{\dot{B}_{2,1}^\fd}}{\frac{c_*}{2}}+ \frac{\|b_m\|_{L^\infty}}{(\frac{c_*}{2})^2}\|\nabla c_m\|_{\dot{B}_{2,1}^\fd}\\
\leq \Cs (\|\nabla b_m\|_{\dot{B}_{2,1}^\fd}+ \|\nabla c_m\|_{\dot{B}_{2,1}^\fd}) = \Cs (\|\nabla (b_m-1)\|_{\dot{B}_{2,1}^\fd}+ \|\nabla (c_m-1)\|_{\dot{B}_{2,1}^\fd}) \\
 \leq 2^m \Cs (\|b-1\|_{\dot{B}_{2,1}^\fd}+ \|c-1\|_{\dot{B}_{2,1}^\fd}).
\label{estimquotientbmcm}
\end{multline}
Therefore, there exists a constant only depending on the bounds of $b$ and $c$ once again denoted by $\Cs$ so that:
\begin{equation}
\Big|(\nabla(v\cdot\nabla q_j)|\frac{b_m}{c_m}\nabla q_j)_{L^2}\Big| \leq \Cs \left(\|\nabla v\|_{L^\infty} + 2^m \Big(\|b-1\|_{\dot{B}_{2,1}^\fd}+ \|c-1\|_{\dot{B}_{2,1}^\fd}\Big) \|v\|_{L^\infty}  \right) \|\nabla q_j\|_{L^2}^2.
\label{estimvgradqj}
\end{equation}
In all the cited works, the right-hand side $(\nabla f_j|\nabla q_j)_{L^2}$ is completely harmless if the coefficients are constant and is estimated  the following way (thanks to the fact that $q_j$ is localized in frequency):
$$
|(\nabla f_j|\nabla q_j)_{L^2}| =|(f_j|\Delta q_j)_{L^2}| \leq C 2^j\|f_j\|_{L^2} \|\nabla q_j\|_{L^2}.
$$
In our case the study of $(\nabla f_j|\frac{b_m}{c_m}\nabla q_j)_{L^2}$ will be more delicate. Indeed $f_j= F_j+F_{m,j}+R_j+\tilde{R}_j$ and as the last two terms are not localized in frequency, as well as $\frac{b_m}{c_m}\nabla q_j$, we will have to be much more careful. We remark that in the variable coefficients cases from \cite{Dbook, Dtruly, Has7, Has8}, as the studies of the density and velocity are decoupled, such a term does not occur.

As $F_j$ and $F_{m,j}$ are localized in frequency we immediately have:
\begin{multline}
\Big|(\nabla F_j +\nabla F_{m,j}|\frac{b_m}{c_m} \nabla q_j)_{L^2}\Big| \leq 2^j (\|F_j\|_{L^2}+ \|F_{m,j}\|_{L^2}) \|\frac{b_m}{c_m}\|_{L^\infty}  \|\nabla q_j\|_{L^2}\\
\leq \Cs 2^j (\|F_j\|_{L^2}+ \|F_{m,j}\|_{L^2}) \|\nabla q_j\|_{L^2}.
\label{fext1}
\end{multline}
The last term $\tilde{R}_j$ is not localized in frequency but is rather easy to estimate: thanks to its definition (we refer to \eqref{deffext}) we have:
\begin{multline}
\tilde{R}_j^l= \d_l \tilde{R}_j = \Sum_{k=1}^d \d_k \left((c_m-1).\d_l u_j^k\right)-\ddj\left((c_m-1).\d_k \d_l u^k\right)\\
+\Sum_{k=1}^d \d_k \left(\d_l(c_m-1) .u_j^k\right)-\ddj\left(\d_l(c_m-1).\d_k u^k\right).
\end{multline}
Using lemma 2 from the appendix of \cite{Dtruly} (see lemma \ref{estimtruly} in the appendix of the present paper), there exists a constant $C>0$ and two nonnegative summable sequences $c_j'(t)$ and $c_j''(t)$ whose summation is 1 such that:
\begin{multline}
\|\nabla \tilde{R}_j\|_{L^2} \leq C c_j' 2^{-j s_1} \|c_m-1\|_{\dot{B}_{2,1}^{\fd+h_1}} \|\nabla u\|_{\dot{B}_{2,1}^{s_1+1-h_1}} \\
+C c_j'' 2^{-j s_2} \|\nabla(c_m-1)\|_{\dot{B}_{2,1}^{\fd+h_2}} \|u\|_{\dot{B}_{2,1}^{s_2+1-h_2}},
\end{multline}
and thanks to the frequency localization of $c_m$,
\begin{multline}
\|\nabla \tilde{R}_j\|_{L^2} \leq C c_j' 2^{-j s_1} 2^{m h_1} \|c-1\|_{\dot{B}_{2,1}^\fd} \|\nabla u\|_{\dot{B}_{2,1}^{s_1+1-h_1}} \\
+C c_j'' 2^{-j s_2} 2^{m(1+h_2)} \|c-1\|_{\dot{B}_{2,1}^\fd} \|u\|_{\dot{B}_{2,1}^{s_2+1-h_2}}.
\end{multline}
With $s_1=s_2=s-1$ and $h_1=h_2=1$, there exists a nonnegative summable sequence $c_j(t)$ whose summation is 1 such that:
$$
\|\nabla \tilde{R}_j\|_{L^2} \leq C c_j 2^{-j (s-1)} \|c-1\|_{\dot{B}_{2,1}^\fd} \left( 2^m  \|\nabla u\|_{\dot{B}_{2,1}^{s-1}} +2^{2m} \|u\|_{\dot{B}_{2,1}^{s-1}} \right).
$$
and then:
\begin{multline}
\Big|(\nabla \tilde{R}_j|\frac{b_m}{c_m} \nabla q_j)_{L^2}\Big| \leq \|\nabla \tilde{R}_j\|_{L^2} \|\frac{b_m}{c_m}\|_{L^\infty} \|\nabla q_j\|_{L^2}\\
\leq \Cs c_j 2^{-j (s-1)} \|c-1\|_{\dot{B}_{2,1}^\fd} \left( 2^m  \|u\|_{\dot{B}_{2,1}^s} +2^{2m} \|u\|_{\dot{B}_{2,1}^{s-1}} \right) \|\nabla q_j\|_{L^2}.
\end{multline}
Finally, thanks to interpolation estimates we obtain
\begin{equation}
\Big|(\nabla \tilde{R}_j|\frac{b_m}{c_m} \nabla q_j)_{L^2}\Big| \leq  c_j 2^{-j (s-1)} \|\nabla q_j\|_{L^2} \left(\frac{1}{2} \|u\|_{\dot{B}_{2,1}^{s+1}} +\Cs 2^{2m} (1+ \|c-1\|_{\dot{B}_{2,1}^\fd} )^2 \|u\|_{\dot{B}_{2,1}^{s-1}} \right).
\label{fext2}
\end{equation}
The real problem comes from
$$
\nabla R_j= \left(v \cdot\nabla^2 q_j-\ddj(v \cdot\nabla^2 q)\right) + \left( \nabla v \cdot\nabla q_j-\ddj(\nabla v \cdot\nabla q)\right).
$$
Indeed, we could use the same estimate as in \cite{Dinv, CH} (see lemma \ref{estimtruly}) for the first two terms, but not for the last two as there are too many derivatives pounding on $v$. We are then forced to make the derivative pound on the second term of the inner product:
$$
(\nabla R_j|\frac{b_m}{c_m} \nabla q_j)_{L^2}= -(R_j|\div(\frac{b_m}{c_m} \nabla q_j))_{L^2} =-(R_j|\frac{b_m}{c_m} \Delta q_j))_{L^2} -(R_j|\nabla q_j \cdot \nabla(\frac{b_m}{c_m} ))_{L^2},
$$
and then, thanks to \eqref{estimquotientbmcm}
\begin{equation}
\Big|(\nabla R_j|\frac{b_m}{c_m} \nabla q_j)_{L^2}\Big| \leq \Cs \left(2^j +2^m(\|b-1\|_{\dot{B}_{2,1}^\fd}+ \|c-1\|_{\dot{B}_{2,1}^\fd}) \right) \|R_j\|_{L^2} \|\nabla q_j\|_{L^2} .
\end{equation}
When $m\leq j$ ($m$ will be more precisely fixed later),
\begin{equation}
\Big|(\nabla R_j|\frac{b_m}{c_m} \nabla q_j)_{L^2}\Big| \leq \Cs \left(1 +(\|b-1\|_{\dot{B}_{2,1}^\fd}+ \|c-1\|_{\dot{B}_{2,1}^\fd}) \right) 2^j \|R_j\|_{L^2} \|\nabla q_j\|_{L^2}.
\label{fext3}
\end{equation}
which gives a good estimate thanks to \eqref{estimfext}. Unfortunately, due to the fact that we only have $q_0 \in \dot{B}_{2,1}^s$ and do not assume $q_0 \in \dot{B}_{2,1}^{s-1}$, this will not work in the case $j \leq m$. In this case we rewrite the equation on $\nabla q_j$ the following way:
$$
\d_t \nabla q_j+\nabla \div(c_m u_j) =\nabla F_j+\nabla F_{m,j}+\nabla \tilde{R}_j -\nabla \ddj (v\cdot\nabla q) ,
$$
and take its innerproduct with $\frac{b_m}{c_m} \nabla q_j$.
\begin{multline}
(\d_t \nabla q_j|\frac{b_m}{c_m}\nabla q_j)_{L^2}+(\nabla \div(c_m u_j)|\frac{b_m}{c_m}\nabla q_j)_{L^2} =\left(\nabla F_j+ \nabla F_{m,j} +\nabla\tilde{R}_j -\nabla \ddj (v\cdot\nabla q)\Big|\frac{b_m}{c_m}\nabla q_j\right)_{L^2}\\
\leq \Cs 2^j (\|F_j\|_{L^2}+ \|F_{m,j}\|_{L^2}) \|\nabla q_j\|_{L^2}  +\left(\nabla \ddj (v\cdot\nabla q)\Big|\frac{b_m}{c_m}\nabla q_j\right)_{L^2}\\
+ c_j 2^{-j (s-1)} \|\nabla q_j\|_{L^2} \left(\frac{1}{2} \|u\|_{\dot{B}_{2,1}^{s+1}} +\Cs 2^{2m} (1+ \|c-1\|_{\dot{B}_{2,1}^\fd} )^2 \|u\|_{\dot{B}_{2,1}^{s-1}} \right).
\label{estimjm}
\end{multline}
Thanks to the Bony decomposition, $v\cdot \nabla q= T_v \nabla q +T_{\nabla q} v + R(v, \nabla q)$ and we can easily show, thanks to the paraproduct and remainder laws for Besov spaces that:
\begin{equation}
\|\ddj\Big(T_{\nabla q} v + R(v, \nabla q)\Big)\|_{L^2} \leq C c_j 2^{-js} \|v\|_{\dot{B}_{2,1}^{\fd+1}} \|q\|_{\dot{B}_{2,1}^s}.
\label{Bony}
\end{equation}
\begin{rem}
\sl{This is here that, due to the paraproduct and remainder conditions on indices, we get the condition $s\in]-\fd, \fd+1]$ in theorem \ref{thestimbc}.}
\end{rem}
But, $T_v \nabla q$ cannot be estimated this way, indeed, the best be can have is:
$$
\|T_v \nabla q\|_{\dot{B}_{2,1}^s} \leq \|v\|_{L^\infty}  \|q\|_{\dot{B}_{2,1}^{s+1}}.
$$
This estimate is useful in the case of (NSK) because the density fluctuation is more regular but not in the non-local case. As we are in the case  $j\leq m$, we will be able to estimate:
$$
T_v \nabla q =\Sum_{l\in \Z} \dot{S}_{l-1} v. \ddl \nabla q.
$$
As $\dot{S}_{l-1} v. \ddl \nabla q$ is localized in frequency in an annulus $2^l \cC'$, using that $q\in \dot{B}_{2,1}^s$ there exists a summable nonnegative sequence whose summation is 1, $(c_l(t))_{l\in \Z}$, such that:
\begin{multline}
2^j \|\ddj(T_v \nabla q)\|_{L^2} \leq 2^j\Sum_{|l-j|\leq N_1} \|\dot{S}_{l-1} v\|_{L^\infty} 2^l \|\ddl q\|_{L^2} \leq C 2^j \Sum_{|l-j|\leq N_1} \|v\|_{L^\infty} 2^l 2^{-ls} c_l(t) \|q\|_{\dot{B}_{2,1}^s}\\
\leq C 2^{-j(s-1)} \Sum_{|l-j|\leq N_1} \|v\|_{L^\infty} 2^l 2^{(j-l)s} c_l(t) \|q\|_{\dot{B}_{2,1}^s}\\
\leq C 2^{m+N_1} 2^{-j(s-1)} \left(\Sum_{|l-j|\leq N_1} 2^{(j-l)s} c_l(t)\right) \|v\|_{L^\infty} \|q\|_{\dot{B}_{2,1}^s}\\
\leq C' 2^m  2^{-j(s-1)} d_j(t) \|v\|_{L^\infty} \|q\|_{\dot{B}_{2,1}^s},
\label{estimvgradqj2}
\end{multline}
because $l\leq j+N_1 \leq m+N_1$ and, thanks to convolution estimates, $(d_j(t))_{j\in \Z}$ is nonnegative, summable with $\|d\|_{l^1} \leq C_s \|c\|_{l^1} \leq C_s$.
Finally we are able to write the desired estimates: collecting \eqref{estimddt}, \eqref{estimdtbmcm} and \eqref{estimvgradqj}, \eqref{fext1}, \eqref{fext2}, \eqref{fext3} (in the case $j \geq m$) and \eqref{fext1}, \eqref{estimjm}, \eqref{Bony}, \eqref{estimvgradqj2} (in the case $j \leq m$) we end up with, for all $j\in \Z$:
\begin{multline}
\frac{1}{2} \frac{d}{dt} (\nabla q_j|\frac{b_m}{c_m}\nabla q_j)_{L^2} +(\nabla \div(c_m u_j)|\frac{b_m}{c_m}\nabla q_j)_{L^2} \leq
\Cs \|\nabla q_j\|_{L^2}^2 \Bigg[ 2^m (\|\d_t b\|_{\dot{B}_{2,1}^{\fd-1}}+ \|\d_t c\|_{\dot{B}_{2,1}^{\fd-1}})\\
+(1+ \|b-1\|_{\dot{B}_{2,1}^\fd}+ \|c-1\|_{\dot{B}_{2,1}^\fd}) (\|\nabla v\|_{\dot{B}_{2,1}^\fd} + 2^m \|v\|_{\dot{B}_{2,1}^\fd})\Bigg]\\
+ \|\nabla q_j\|_{L^2} \Bigg[ c_j 2^{-j(s-1)} \|u\|_{\dot{B}_{2,1}^{s+1}} +\Cs 2^j (\|F_j\|_{L^2} +\|F_{m,j}\|_{L^2})\\
+c_j 2^{-j(s-1)} \Cs 2^{2m} (1+\|c-1\|_{\dot{B}_{2,1}^\fd})^2 \|u\|_{\dot{B}_{2,1}^{s-1}}\\
+ c_j 2^{-j(s-1)} \Cs (1+ \|b-1\|_{\dot{B}_{2,1}^\fd}+ \|c-1\|_{\dot{B}_{2,1}^\fd}) (\|\nabla v\|_{\dot{B}_{2,1}^\fd} + 2^m \|v\|_{\dot{B}_{2,1}^\fd}) \|q\|_{\dot{B}_{2,1}^s}\Bigg].
\label{estimgradqj2}
\end{multline}
Let us recall that the interest of computing \eqref{estimujgradqj}$+\nu$\eqref{estimgradqj2} is that, as in \cite{CH}, it allows to neutralize \eqref{difftermes}, that is the terms consuming too many derivatives (estimated by $\Cs 2^{2j}\|u_j\|_{L^2}\|\nabla q_j\|_{L^2}$), and in the remaining terms at least one derivative will pound on $b_m-1$ or $c_m-1$ making these terms harmless because there is $2^m 2^j$ or $2^{2m}$ instead of $2^{2j}$. Let us explain this in details: as $\nu= \lambda +2\mu$, we have (from \eqref{estimujgradqj}$+\nu$\eqref{estimgradqj2})
\begin{multline}
B\overset{def}{=}-\mu (\div(b_m.\nabla u_j)| \nabla q_j)_{L^2}-(\lambda+\mu) (\nabla(b_m \div u_j)| \nabla q_j)_{L^2}
+\nu (\nabla \div(c_m u_j)|\frac{b_m}{c_m}\nabla q_j)_{L^2}\\
=\mu\left[(\nabla \div(c_m u_j)|\frac{b_m}{c_m}\nabla q_j)_{L^2} -(\div(b_m.\nabla u_j)| \nabla q_j)_{L^2} \right]\\
+(\lambda+\mu) \left[(\nabla \div(c_m u_j)|\frac{b_m}{c_m}\nabla q_j)_{L^2} -(\nabla(b_m \div u_j)| \nabla q_j)_{L^2}\right] = \mu B_1 +(\lambda+\mu) B_2.
\end{multline}
Let us now estimate separatedly $B_1$ and $B_2$:
\begin{multline}
B_2= \left(\frac{b_m}{c_m} (\nabla c_m\div u_j +c_m \nabla\div u_j  +\nabla u_j.\nabla c_m +u_j.\nabla^2 c_m)\Big|\nabla q_j\right)_{L^2}\\ -(b_m\nabla \div u_j +\nabla b_m.\div u_j | \nabla q_j)_{L^2}\\
=\left(\frac{b_m}{c_m} \nabla c_m\div u_j +\frac{b_m}{c_m}\nabla c_m .\nabla u_j +\frac{b_m}{c_m}.\nabla^2 c_m .u_j -\nabla b_m.\div u_j  \Big|\nabla q_j\right)_{L^2}.\\
\end{multline}
We emphasize that $(b_m\nabla \div u_j)| \nabla q_j)_{L^2}$, which was mentionned before as an obstruction term, has disappeared, and all that remain are harmless terms where at most two derivatives are applied to $u_j$ or $q_j$. Thanks to \eqref{estimquotientbmcm}, and the fact that $\nabla^2 c_m =\nabla^2 (c_m-1)= \nabla^2 \dot{S}_m (c-1)$ we obtain
$$
|B_2| \leq \Cs (\|b-1\|_{\dot{B}_{2,1}^\fd} +\|c-1\|_{\dot{B}_{2,1}^\fd}) \left( 2^m  2^j \|u_j\|_{L^2}+ 2^{2m} \|u_j\|_{L^2}\right) \|\nabla q_j\|_{L^2}.
$$
The other term requires a little more attention:
\begin{multline}
B_1= -\left(\div(c_m u_j)\Big|\div(\frac{b_m}{c_m}\nabla q_j)\right)_{L^2} -(b_m\Delta u_j+ \nabla u_j. \nabla b_m)| \nabla q_j)_{L^2}\\
=-\left(c_m\div u_j+ u_j.\nabla c_m\Big|\frac{b_m}{c_m} \Delta q_j +\nabla q_j.\nabla (\frac{b_m}{c_m})\right)_{L^2} -(b_m\Delta u_j+ \nabla u_j. \nabla b_m)| \nabla q_j)_{L^2}\\
=B_{11} +B_{12},
\end{multline}
with
$$
B_{11}=- (b_m\div u_j|\Delta q_j)_{L^2} -(b_m\Delta u_j|\nabla q_j)_{L^2},
$$
and
\begin{multline}
B_{12}=-\left(u_j.\nabla c_m\Big|\frac{b_m}{c_m} \Delta q_j +\nabla q_j.\nabla (\frac{b_m}{c_m})\right)_{L^2} -\left(c_m\div u_j \Big|\nabla q_j.\nabla (\frac{b_m}{c_m})\right)_{L^2}\\
-(\nabla u_j. \nabla b_m)| \nabla q_j)_{L^2}.
\end{multline}
In each term of $B_{12}$ at most two derivatives act on $u_j$ and $q_j$ so that we can estimate $B_{12}$ with the same arguments used for $B_2$, and for the other term we have:
\begin{multline}
B_{11}=- \left(\Delta(b_m\div u_j)\Big| q_j\right)_{L^2} -(b_m\Delta u_j|\nabla q_j)_{L^2}\\
=- (\Delta b_m.\div u_j + 2\nabla b_m.\nabla \div u_j +b_m. \Delta \div u_j  | q_j )_{L^2} +(b_m.\div \Delta u_j +\Delta u_j.\nabla b_m| q_j)_{L^2}\\
=- (\Delta b_m.\div u_j + 2\nabla b_m.\nabla \div u_j +\Delta u_j.\nabla b_m | q_j )_{L^2}.
\end{multline}
As expected, the dangerous terms are neutralized and then, as before, we can easily estimate $B_1$ and we can finally write:
\begin{multline}
|B| \leq \Cs \overline{\nu}(\|b-1\|_{\dot{B}_{2,1}^\fd} +\|c-1\|_{\dot{B}_{2,1}^\fd})\\
\times\left( 2^m  2^j \|u_j\|_{L^2} \|\nabla q_j\|_{L^2}+ 2^{2m} \Big(1+\|c-1\|_{\dot{B}_{2,1}^\fd}\Big) \|u_j\|_{L^2} \|\nabla q_j\|_{L^2}\right),
\end{multline}
and then, using the fact that $2ab\leq a^2 +b^2$, we finally end up with:
\begin{multline}
|B| \leq \Cs 2^{2j} \|u_j\|_{L^2}^2 +\Cs \overline{\nu}^2 2^{2m} (1+\|b-1\|_{\dot{B}_{2,1}^\fd} +\|c-1\|_{\dot{B}_{2,1}^\fd})^2 \left(\|u_j\|_{L^2} \|\nabla q_j\|_{L^2} +\|\nabla q_j\|_{L^2}^2\right).
\label{estimB}
\end{multline}
Thanks to \eqref{estimB}, computing \eqref{estimujgradqj}$+\nu$\eqref{estimgradqj2} allows us to obtain (using once more \eqref{estimfext}, and estimating $2^j \|\tilde{R}_j\|_{L^2}$ and $\|\tilde{S}_j\|_{L^2}$ like in \eqref{fext2}), denoting by $\Ds$ a constant only depending on the bounds of $b,c$ and on the physical parameters $\mu$, $\lambda$, $\kappa$:
\begin{multline}
\frac{d}{dt} \left( (u_j|\nabla q_j)_{L^2} +\frac{\nu}{2} (\nabla q_j|\frac{b_m}{c_m}\nabla q_j)_{L^2} \right) +\kappa (\nabla q_j|\frac{\nabla q_j-\phie *\nabla q_j}{\ee^2})_{L^2} \leq \Ds 2^{2j}\|u_j\|_{L^2}^2\\
+ \Ds (\|u_j\|_{L^2}^2 +\|\nabla q_j\|_{L^2}^2) \Bigg[ 2^m (\|\d_t b\|_{\dot{B}_{2,1}^{\fd-1}}+ \|\d_t c\|_{\dot{B}_{2,1}^{\fd-1}})\\
+(1+ \|b-1\|_{\dot{B}_{2,1}^\fd}+ \|c-1\|_{\dot{B}_{2,1}^\fd})^2 (\|\nabla v\|_{\dot{B}_{2,1}^\fd} + 2^m \|v\|_{\dot{B}_{2,1}^\fd} +2^{2m})\Bigg]\\
+\Ds (\|u_j\|_{L^2} +\|\nabla q_j\|_{L^2}) \Bigg[ c_j 2^{-j(s-1)} \|u\|_{\dot{B}_{2,1}^{s+1}} +(\|G_j\|_{L^2} +\|G_{m,j}\|_{L^2} +2^j \|F_j\|_{L^2} +2^j \|F_{m,j}\|_{L^2})\\
+c_j 2^{-j(s-1)} (1 +\|b-1\|_{\dot{B}_{2,1}^\fd} +\|c-1\|_{\dot{B}_{2,1}^\fd})^2 (\|\nabla v\|_{\dot{B}_{2,1}^\fd} +2^{2m}) \|u\|_{\dot{B}_{2,1}^{s-1}}\\
+ c_j 2^{-j(s-1)} (1+ \|b-1\|_{\dot{B}_{2,1}^\fd}+ \|c-1\|_{\dot{B}_{2,1}^\fd}) (\|\nabla v\|_{\dot{B}_{2,1}^\fd} + 2^m \|v\|_{\dot{B}_{2,1}^\fd}) \|q\|_{\dot{B}_{2,1}^s}\Bigg].
\label{estimujqj2}
\end{multline}
Let us introduce, as stated in theorem \ref{thestimbc}, the following quantity
\begin{equation}
h_j(t)^2 =(u_j|c_m u_j)_{L^2} +\kappa (q_j|\frac{q_j-\phie *q_j}{\ee^2})_{L^2} +\eta \left(2 (u_j|\nabla q_j)_{L^2} +\nu (\nabla q_j|\frac{b_m}{c_m} \nabla q_j)_{L^2}  \right),
\end{equation}
for some $\eta<1$ that we will precise in the following. First we will write conditions on $\eta$ so that $h_j^2 \sim \|u_j\|_{L^2}^2 +\|\nabla q_j\|_{L^2}^2$. Thanks to \eqref{bornesbc},
$$
\nu \frac{b_*}{2c^*+c_*} \|\nabla q_j\|_{L^2}^2 \leq \nu (\nabla q_j| \frac{b_m}{c_m} \nabla q_j)_{L^2} \leq \nu \frac{2b^*+b_*}{c_*} \|\nabla q_j\|_{L^2}^2.
$$
Using the fact that $2ab\leq a^2 +b^2$, we get
$$
2|(u_j|\nabla q_j)_{L^2}| \leq \nu \frac{b_*}{2(2c^*+c_*)} \|\nabla q_j\|_{L^2}^2 + \frac{2(2c^*+c_*)}{\nu b_*} \|u_j\|_{L^2}^2
$$
so that if $\eta$ satisfies
\begin{equation}
\eta \leq \min(1,\frac{\nu b_*c_*}{8 (2c^*+c_*)}),
\label{condeta1}
\end{equation}
there exists a constant $\Ds$ such that we have:
\begin{multline}
\frac{c_*}{4} \|u_j\|_{L^2}^2 +\kappa (q_j|\frac{q_j-\phie *q_j}{\ee^2})_{L^2} +\frac{\nu \eta b_*}{2(2c^*+c_*)} \|\nabla q_j\|_{L^2}^2 \leq h_j(t)^2\\
 \leq \Ds (\|u_j\|_{L^2}^2 +\|\nabla q_j\|_{L^2}^2) +\kappa (q_j|\frac{q_j-\phie *q_j}{\ee^2})_{L^2}.
 \label{equivhj}
\end{multline}
then, writing \eqref{estimujqj1}$+2\eta$\eqref{estimujqj2}, and using the fact that (in \eqref{estimujqj1}):
\begin{multline}
C \kappa \|v\|_{L^\infty} \|\nabla q_j\|_{L^2} 2^j \sqrt{(q_j| \frac{q_j -\phie *q_j}{\ee^2})_{L^2}}\\
\leq \frac{\kappa \eta}{2} (\nabla q_j|\frac{\nabla q_j-\phie *\nabla q_j}{\ee^2})_{L^2} +\frac{C^2 \kappa}{2 \eta} \|v\|_{L^\infty}^2 \|\nabla q_j\|_{L^2}^2,
\end{multline}
if $\eta$ satisfies
\begin{equation}
\eta \leq \frac{\nu_0 b_* c_*}{16 \Ds},
\label{condeta2}
\end{equation}
then for all $t\in[0,T]$:
\begin{multline}
\frac{1}{2} \frac{d}{dt} h_j(t)^2 +\nu_0 \frac{b_* c_*}{16} 2^{2j} \|u_j\|_{L^2}^2 +\frac{\kappa \eta}{2}(\nabla q_j|\frac{\nabla q_j-\phie *\nabla q_j}{\ee^2})_{L^2}\\
\leq \Ds h_j(t)^2\Bigg[ 2^m (\|\d_t b\|_{\dot{B}_{2,1}^{\fd-1}}+ \|\d_t c\|_{\dot{B}_{2,1}^{\fd-1}})\\
+(1+ \|b-1\|_{\dot{B}_{2,1}^\fd}+ \|c-1\|_{\dot{B}_{2,1}^\fd})^2 (\|\nabla v\|_{\dot{B}_{2,1}^\fd} + 2^m \|v\|_{\dot{B}_{2,1}^\fd} +2^{2m} +\|v\|_{\dot{B}_{2,1}^\fd}^2 )\Bigg]\\
+\Ds h_j(t) \Bigg[ c_j 2^{-j(s-1)} \eta \|u\|_{\dot{B}_{2,1}^{s+1}} +(\|G_j\|_{L^2} +\|G_{m,j}\|_{L^2} +2^j \|F_j\|_{L^2} +2^j \|F_{m,j}\|_{L^2})\\
+c_j 2^{-j(s-1)} (1 +\|b-1\|_{\dot{B}_{2,1}^\fd} +\|c-1\|_{\dot{B}_{2,1}^\fd})^2 (\|\nabla v\|_{\dot{B}_{2,1}^\fd} +2^{2m}) \|u\|_{\dot{B}_{2,1}^{s-1}}\\
+ c_j 2^{-j(s-1)} (1+ \|b-1\|_{\dot{B}_{2,1}^\fd}+ \|c-1\|_{\dot{B}_{2,1}^\fd}) (\|\nabla v\|_{\dot{B}_{2,1}^\fd} + 2^m \|v\|_{\dot{B}_{2,1}^\fd}) \|q\|_{\dot{B}_{2,1}^s}\Bigg].
\end{multline}
If we introduce
$$
\gamma_*= \min(\frac{\kappa \eta}{4 \Ds}, \frac{\nu_0 b_* c_*}{16 \Ds}, \frac{\eta}{4}),
$$
then for all $j\in \Z$,
$$
\nu_0 \frac{b_* c_*}{16} 2^{2j} \|u_j\|_{L^2}^2 +\frac{\kappa \eta}{2}(\nabla q_j|\frac{\nabla q_j-\phie *\nabla q_j}{\ee^2})_{L^2}
\geq \gamma_* \min(\frac{1}{\ee^2}, 2^{2j}) h_j(t)^2,
$$
and we obtain, integrating from $0$ to $t$:
\begin{multline}
h_j(t) +\gamma_* \min(\frac{1}{\ee^2}, 2^{2j}) \int_0^t h_j(\tau) d\tau \leq h_j(0)\\
+ \Ds \int_0^t h_j(\tau)\Bigg[ 2^m (\|\d_t b\|_{\dot{B}_{2,1}^{\fd-1}}+ \|\d_t c\|_{\dot{B}_{2,1}^{\fd-1}})\\
+(1+ \|b-1\|_{\dot{B}_{2,1}^\fd}+ \|c-1\|_{\dot{B}_{2,1}^\fd})^2 (\|\nabla v\|_{\dot{B}_{2,1}^\fd} + 2^m \|v\|_{\dot{B}_{2,1}^\fd} +2^{2m} +\|v\|_{\dot{B}_{2,1}^\fd}^2 )\Bigg] d\tau\\
+\Ds \int_0^t \Bigg[ c_j 2^{-j(s-1)} \eta \|u\|_{\dot{B}_{2,1}^{s+1}} +(\|G_j\|_{L^2} +\|G_{m,j}\|_{L^2} +2^j \|F_j\|_{L^2} +2^j \|F_{m,j}\|_{L^2})\\
+c_j 2^{-j(s-1)} (1 +\|b-1\|_{\dot{B}_{2,1}^\fd} +\|c-1\|_{\dot{B}_{2,1}^\fd})^2 (\|\nabla v\|_{\dot{B}_{2,1}^\fd} +2^{2m}) \|u\|_{\dot{B}_{2,1}^{s-1}}\\
+ c_j 2^{-j(s-1)} (1+ \|b-1\|_{\dot{B}_{2,1}^\fd}+ \|c-1\|_{\dot{B}_{2,1}^\fd}) (\|\nabla v\|_{\dot{B}_{2,1}^\fd} + 2^m \|v\|_{\dot{B}_{2,1}^\fd}) \|q\|_{\dot{B}_{2,1}^s}\Bigg]d\tau.
\end{multline}
In particular, we obtained a bound for the integral (see remark \ref{rqintegraleepsilon}) :
$$
\min(\frac{1}{\ee^2}, 2^{2j}) \int_0^t \|\nabla q_j(\tau)\|_{L^2} d\tau,
$$
that is, we can now obtain the parabolic regularization on the velocity exactly as in \cite{Dinv, CH} by taking the innerproduct of the velocity equation from system \eqref{systloc} by $u_j$, finally if we set, as stated in theorem \ref{thestimbc}, $g_j(t) =\|u_j(t)\|_{L^2} +h_j(t)$, then we have \eqref{equivgs}, and multiplying by $2^{j(s-1)}$ and summing over $j\in \Z$, we obtain:
\begin{multline}
g^s(q,u)(t) +\nu_0 \frac{b_*}{2} 2^{2j} \int_0^t \|u\|_{\dot{B}_{2,1}^{s+1}} d\tau +\gamma_* \int_0^t \|\frac{\phie*q-q}{\ee^2}\|_{\dot{B}_{2,1}^s}  d\tau\\
\leq g^s(q_0,u_0)(0) +\Ds \int_0^t g^s(q,u)(\tau)\Bigg[ 2^m (\|\d_t b\|_{\dot{B}_{2,1}^{\fd-1}}+ \|\d_t c\|_{\dot{B}_{2,1}^{\fd-1}})\\
+(1+ \|b-1\|_{\dot{B}_{2,1}^\fd}+ \|c-1\|_{\dot{B}_{2,1}^\fd})^2 (\|\nabla v\|_{\dot{B}_{2,1}^\fd} + 2^m \|v\|_{\dot{B}_{2,1}^\fd} +2^{2m} +\|v\|_{\dot{B}_{2,1}^\fd}^2 )\Bigg] d\tau\\
+\Ds \int_0^t \Bigg[ \eta \|u\|_{\dot{B}_{2,1}^{s+1}} +(\|G\|_{\dot{B}_{2,1}^{s-1}} +\|G_m\|_{\dot{B}_{2,1}^{s-1}} +\|F\|_{\dot{B}_{2,1}^s} +\|F_m\|_{\dot{B}_{2,1}^s})\Bigg] d\tau.
\end{multline}
Using \eqref{estimfext} and finally fixing $\eta$ such that
\begin{equation}
\eta \leq \frac{\nu_0 b_*}{4\Ds},
\label{condeta3}
\end{equation}
we end up with \eqref{estimgs}.
\begin{rem}
\sl{Recollecting the various conditions on $\eta$, \eqref{condeta1}, \eqref{condeta2} and \eqref{condeta3}, we need that:
$$
0< \eta \leq \min(1,\frac{\nu b_*c_*}{8 (2c^*+c_*)}, \frac{\nu_0 b_* c_*}{16 \Ds}, \frac{\nu_0 b_*}{4\Ds}),
$$
and then we defined
$$
\gamma_*= \min(\frac{\kappa \eta}{4 \Ds}, \frac{\nu_0 b_* c_*}{16 \Ds}, \frac{\eta}{4})>0,
$$
where $\Ds$ only depends on the bounds of $b$ and $c$, the physical parameters (viscosities, capillarity, dimension) and $s$.
}
\end{rem}

\begin{rem}
\sl{As explained, many points of the previous proof are much easier for the Korteweg model.}
\end{rem}

\subsection{Existence}

This part is classical so we will leave some details to the reader (the proofs are very similar to those in \cite{Dheatlocal, Dheatglobal, Dtruly, Dbook, CH, Has7, Has8}). As in \cite{Dtruly}, we will first prove the existence for more regular initial data (we refer to \cite{Dheatlocal, Dheatglobal}), then use it for regularized initial data, and finally prove the convergence to a solution of our problem.

The special feature of the present models is that, due to the capillary term, we cannot afford to study separatedly the density fluctuation and the velocity, we need to recouple them via the previous a priori estimate. Again we present the proof in the non-local case (it is easier for the Korteweg system).

Let us first assume that the initial data satisfy $q_0 \in \dot{B}_{2,1}^\fd \cap \dot{B}_{2,1}^{\fd+1}$ and $u_0 \in \dot{B}_{2,1}^{\fd-1} \cap \dot{B}_{2,1}^\fd$. We emphasize that this will not be necessary for system (NSK). We introduce $(q_L, u_L)$ the unique global solution of system:
\begin{equation}
\begin{cases}
\begin{aligned}
&\d_t q_L +\div u_L= 0,\\
&\d_t u_L -\cA u_L-\kappa \frac{\phie*\nabla q_L-\nabla q_L}{\ee^2}= 0,\\
&(q_L, u_L)_{|t=0}=(q_0, u_0).
\end{aligned}
\end{cases}
\label{systL}
\end{equation}
We can easily adapt the linear estimates from \cite{CHVP} in the case where we do not assume $P'(1)>0$ (that is taking $p=0$ in system $(LR_\ee)$ from \cite{CH} or \cite{CHVP}), we obtain the following global in time estimates (let us recall that we use the estimates from \cite{DD} for $(NSK)$, or adapt those from \cite{Corder} for the order parameter model): there exists a constant $C>0$ such that for all $t$, and for $s\in\{\fd, \fd+1\}$,
\begin{multline}
\nu \|q_L\|_{\tilde{L}^\infty_t\dot{B}_{2,1}^s} +\nu^2\|q_L\|_{L^1_t\dot{B}_{1/\ee}^{s+2,s}} +\|u_L\|_{\tilde{L}^\infty_t\dot{B}_{2,1}^{s-1}} +\nu_0\|u_L\|_{L^1_t\dot{B}_{2,1}^{s+1}}\\
\leq C(\nu \|q_0\|_{\dot{B}_{2,1}^s} +\|u_0\|_{\dot{B}_{2,1}^{s-1}}).
\end{multline}
As this system is linear, for all $m\in \Z$, thanks to this estimates we also have for all $t$,
\begin{multline}
\nu \|(I_d-\dot{S}_m)q_L\|_{\tilde{L}^\infty_t\dot{B}_{2,1}^s} +\nu^2\|(I_d-\dot{S}_m)q_L\|_{L^1_t\dot{B}_{1/\ee}^{s+2,s}} +\|(I_d-\dot{S}_m)u_L\|_{\tilde{L}^\infty_t\dot{B}_{2,1}^{s-1}}\\
+\nu_0\|(I_d-\dot{S}_m)u_L\|_{L^1_t\dot{B}_{2,1}^{s+1}} \leq C(\nu \|(I_d-\dot{S}_m)q_0\|_{\dot{B}_{2,1}^s} +\|(I_d-\dot{S}_m)u_0\|_{\dot{B}_{2,1}^{s-1}}).
\end{multline}
\begin{rem}
\sl{Of course, we could also use our a priori estimates in the case $v=0$, $b=c=1$ on the time interval $[0,T]$ provided by theorem \ref{thexistR}, the interest of the previous estimates being the precision of the coefficients.}
\end{rem}
As in \cite{Dtruly} we define the following iterative scheme: for $n\geq 0$, $(q_{n+1}, u_{n+1})$ is the unique global solution of the following  linear system:
\begin{equation}
\begin{cases}
\begin{aligned}
&\d_t q_{n+1}+ u_n\cdot\nabla q_{n+1}+ (1+q_n)\div u_{n+1}=0,\\
&\d_t u_{n+1}+ u_n\cdot\nabla u_{n+1} -\frac{1}{1+q_n}\cA u_{n+1}-\kappa \frac{\phie*\nabla q_{n+1}-\nabla q_{n+1}}{\ee^2}=-\nabla(H(1+q_n)-H(1)),\\
&(q_{n+1}, u_{n+1})_{|t=0}=(q_0, u_0),
\end{aligned}
\end{cases}
\label{schema}
\end{equation}
where $H$ is the real-valued function defined on $\R$ by $H'(x)=P'(x)/x$. We also define the difference $(\overline{q}_{n+1}, \overline{u}_{n+1})=(q_{n+1}-q_L, u_{n+1}-u_L)$ satisfies the following system:
\begin{equation}
\begin{cases}
\begin{aligned}
&\d_t \overline{q}_{n+1}+ u_n\cdot\nabla \overline{q}_{n+1}+ (1+q_n)\div \overline{u}_{n+1} \quad =-u_n \cdot \nabla q_L -q_n \div u_L =\overline{F}_n,\\
&\d_t \overline{u}_{n+1}+ u_n\cdot\nabla \overline{u}_{n+1} -\frac{1}{1+q_n}\cA \overline{u}_{n+1}-\kappa \frac{\phie*\nabla \overline{q}_{n+1}-\nabla \overline{q}_{n+1}}{\ee^2}\\
& \hspace{5cm} =-\nabla(H(1+q_n)-H(1))-u_n \cdot \nabla u_L -\frac{q_n}{1+q_n} \cA u_L =\overline{G}_n,\\
&(\overline{q}_{n+1}, \overline{u}_{n+1})_{|t=0}=(0, 0),
\end{aligned}
\end{cases}
\label{schemadiff}
\end{equation}
Let us denote by $J$ the real function $J(x)=1/(1+x)$ and assume that $\cu \leq 1+q_0 \leq \co$ (then $\frac{1}{\co} \leq J(q_0) \leq \frac{1}{\cu}$). Let us define $\Fs, \gamma_*$ and $\eta$ as provided by theorem \ref{thestimbc} applied for functions $b,c$ having the following bounds:
$$
c_*=\frac{\cu}{2} \leq c \leq c^*=\co +\frac{\cu}{2}, \quad \mbox{and} \quad b_*=\frac{1}{c^*} \leq b \leq b^*=\frac{1}{c_*}.
$$
For $n \in \N^*$, $t>0$, $m\in \Z$, $0<2\beta \leq \min(b_*,c_*)$ and two positive constants $C, C_0$, we define the following proposition:
\begin{equation}
\mathcal{P}(n)=
\begin{cases}
(1) \quad g^s(q_n, u_n)(T) +\frac{\nu_0 b_*}{4} \|u_n\|_{L_T^1 \dot{B}_{2,1}^{s+1}} +\gamma^* \|q_n\|_{L^1_T\dot{B}_{1/\ee}^{s+2,s}} \leq 2C_0, \quad \mbox{ for } s\in\{\fd,\fd+1\},\\
(2) \quad g^\fd(\overline{q}_n, \overline{u}_n)(T) +\frac{\nu_0 b_*}{4} \|\overline{u}_n\|_{L_T^1 \dot{B}_{2,1}^{\fd+1}} +\gamma^* \|\overline{q}_n\|_{L^1_T\dot{B}_{1/\ee}^{\fd+2,\fd}} \leq \beta,\\
(3) \quad \frac{\nu_0 b_*}{4} \|u_n\|_{L_T^1 \dot{B}_{2,1}^{\fd+1}} +\gamma^* \|q_n\|_{L^1_T\dot{B}_{1/\ee}^{\fd+2,\fd}} \leq 2\beta,\\
(4) \quad \|\d_t q_n\|_{L_T^2 \dot{B}_{2,1}^{\fd-1}} \leq 2C \Fs (1+2C_0)^2,\\
(5) \quad \|\d_t J(q_n)\|_{L_T^2 \dot{B}_{2,1}^{\fd-1}} \leq 2C \Fs (1+2C_0)^4,\\
(6) \quad \|(I_d-\dot{S}_m)q_n\|_{\tilde{L}^\infty_T\dot{B}_{2,1}^\fd} \leq 2 \beta,\\
(7) \quad \|q_n-q_0\|_{\tilde{L}^\infty_T\dot{B}_{2,1}^\fd} \leq 2\beta \quad \mbox{ and } c_* \leq 1+q_n \leq c^*, \quad b_* \leq J(q_n) \leq b^*,\\
(8) \quad \|J(q_n)-1\|_{\tilde{L}^\infty_T\dot{B}_{2,1}^\fd} \leq 2C \Fs C_0,\\
(9) \quad \|J(q_n)-J(q_0)\|_{\tilde{L}^\infty_T\dot{B}_{2,1}^\fd} \leq 2\beta,\\
(10) \quad \|(I_d-\dot{S}_m)\left(J(q_n)-1\right)\|_{\tilde{L}^\infty_T\dot{B}_{2,1}^\fd} \leq 2 \beta
\end{cases}
\end{equation} 
Let us state the main ingredient for the proof of the existence theorem:
\begin{lem}
\sl{Let $C_0=\max(g^\fd(q_0, u_0)(0), g^{\fd+1}(q_0, u_0)(0), \|J(q_0)-1\|_{\dot{B}_{2,1}^\fd})$. Let $\Fs$, $\gamma_*$, $b_*$, $b^*$, $c_*$, $c^*$ as before, $0<\beta< \min(b_*,c_*, \frac{\nu_0 b_*}{32 \Fs})$ small, and assume that $m\in \Z$ is large enough so that:
\begin{equation}
\|(I_d-\dot{S}_m)q_0\|_{\tilde{L}^\infty_t\dot{B}_{2,1}^\fd} +\|(I_d-\dot{S}_m)\left(J(q_0)-1\right)\|_{\tilde{L}^\infty_t\dot{B}_{2,1}^\fd} \leq \beta.
\label{condm}
\end{equation}
If $\beta$ then $T$ are small enough so that:
\begin{equation}
\begin{cases}
e^{2C\Fs\left(1+2(C+1)C_0\right)^4 (1+2 C_0)^4 \beta} \leq \frac{5}{4},\\
\frac{5}{4} (1+\Fs T) e^{2C\Fs\left(1+2(C+1)C_0\right)^4 \left(2^{2m+1}T+ 2^m T^{\frac{1}{2}}\right)} \leq 2,
\end{cases}
\label{condbetaT}
\end{equation}
and
\begin{equation}
C_0 \Fs \left( \|q_L\|_{L_T^2 \dot{B}_{2,1}^{\fd+1}} +\|u_L\|_{L_T^1 \dot{B}_{2,1}^{\fd+1}} +\|u_L\|_{L_T^2 \dot{B}_{2,1}^{\fd}}\right) \leq \beta,
\label{condbetaT2}
\end{equation}
then $\mathcal{P}(n)$ is true for all $n\geq 1$.
}
\end{lem}
\begin{rem}
\sl{We emphasize that as for all $n$, $1+q_n$ and $J(q_n)$ have the same bounds, the constants $\Fs$ and $\gamma_*$ provided by theorem \ref{thestimbc} are the same for all $n$.}
\end{rem}
\begin{rem}
\sl{ Note that more precisely, in (1), in the case $s=\fd$, the constant $C_0$ only depends on norms with regularity index $\fd$.}
\end{rem}
\textbf{Proof:} We will prove this result by induction. Assuming that the proposition is true for $n\in \N^*$, we first apply theorem \ref{thestimbc} to system \eqref{schema}. We have $F=0$ and thanks to composition estimates in Besov spaces (see the appendix), for $s\in\{\fd, \fd+1\}$
$$
\|G\|_{\dot{B}_{2,1}^{s-1}} \leq C(\|q_n\|_{L^\infty})\|q_n\|_{\dot{B}_{2,1}^s} \leq C(C.\|q_0\|_{L^\infty}) \Fs 2C_0,
$$
so that for all $t\in [0,T]$ (thanks to $(1,4,5,6,10)$),
\begin{multline}
g^s(q_{n+1}, u_{n+1})(t) +\frac{\nu_0 b_*}{4} \|u_{n+1}\|_{L_t^1 \dot{B}_{2,1}^{s+1}} +\gamma^* \|q_{n+1}\|_{L^1_t\dot{B}_{\frac{1}{\ee}}^{s+2,s}} \leq g^s(q_0, u_0)(0)\\
+\Fs C_0 t +3\beta \Fs \|u_{n+1}\|_{L_t^1 \dot{B}_{2,1}^{s+1}}+\Fs \int_0^t g^s(q_{n+1}, u_{n+1})(\tau) \Bigg(4C(1+2C_0)^4 2^m\\ +(1+2(C+1)C_0)^2 \Big[2^m\|u_n\|_{\dot{B}_{2,1}^\fd} +\|u_n\|_{\dot{B}_{2,1}^\fd}^2 +\|u_n\|_{\dot{B}_{2,1}^{\fd+1}} +2^{2m}\Big]\Bigg) d\tau,\\
\end{multline}
then using that $3\beta \Fs \leq \frac{\nu_0 b_*}{8}$ and the Gronwall lemma, we obtain:
\begin{multline}
g^s(q_{n+1}, u_{n+1})(t) +\frac{\nu_0 b_*}{4} \|u_{n+1}\|_{L_t^1 \dot{B}_{2,1}^{s+1}} +\gamma^* \|q_{n+1}\|_{L^1_t\dot{B}_{\frac{1}{\ee}}^{s+2,s}} \leq \left(g^s(q_0, u_0)(0) +\Fs C_0 t\right)\\
\times e^{\Fs 2C(1+2(C+1)C_0)^4 \int_0^t \Big[ 2^{m+1} +2^{2m} +2^m\|u_n\|_{\dot{B}_{2,1}^\fd} +\|u_n\|_{\dot{B}_{2,1}^\fd}^2 +\|u_n\|_{\dot{B}_{2,1}^{\fd+1}} \Big] d\tau}.
\end{multline}
Using $(1)$ for $s=\fd$ and $(3)$ from $\mathcal{P}(n)$, we obtain (denoting again by $\Fs$ the constants):
\begin{multline}
g^s(q_{n+1}, u_{n+1})(t) +\frac{\nu_0 b_*}{4} \|u_{n+1}\|_{L_t^1 \dot{B}_{2,1}^{s+1}} +\gamma^* \|q_{n+1}\|_{L^1_t\dot{B}_{\frac{1}{\ee}}^{s+2,s}}\\
\leq C_0(1+\Fs t) e^{\Fs 2C(1+2(C+1)C_0)^4 \Big[ 2^{2m+1}t+2^m t^{\frac{1}{2}} \beta^{\frac{1}{2}} +\beta(1+2C_0)\Big] d\tau}.
\end{multline}
When $\beta>0$ and $T>0$ are small enough to fulfill condition \eqref{condbetaT} ($m$ being fixed more precisely later) we obtain point $(1)$ of $\mathcal{P}(n+1)$.
Applying theorem \ref{thestimbc} to system \eqref{schemadiff} and doing the same leads to (thanks to \eqref{condbetaT}):
\begin{multline}
g^\fd(\overline{q}_{n+1}, \overline{u}_{n+1})(t) +\frac{\nu_0 b_*}{4} \|\overline{u}_{n+1}\|_{L_t^1 \dot{B}_{2,1}^{\fd+1}} +\gamma^* \|\overline{q}_{n+1}\|_{L^1_t\dot{B}_{\frac{1}{\ee}}^{\fd+2,\fd}}\\
\leq e^{\Fs 2C(1+2(C+1)C_0)^4 \Big[ 2^{2m+1}t+2^m t^{\frac{1}{2}} \beta^{\frac{1}{2}} +\beta(1+2C_0)\Big] d\tau} \Fs \int_0^t \left(\|\overline{F}_n\|_{\dot{B}_{2,1}^\fd} +\|\overline{G}_n\|_{\dot{B}_{2,1}^{\fd-1}}  \right) d\tau\\
\leq 2\Fs  \int_0^t \left(\|\overline{F}_n\|_{\dot{B}_{2,1}^\fd} +\|\overline{G}_n\|_{\dot{B}_{2,1}^{\fd-1}}  \right) d\tau.\\
\end{multline}
Estimating the forcing terms:
\begin{multline}
\int_0^t \left(\|\overline{F}_n\|_{\dot{B}_{2,1}^\fd} +\|\overline{G}_n\|_{\dot{B}_{2,1}^{\fd-1}}  \right) d\tau \leq \|u_n\|_{L_t^2 \dot{B}_{2,1}^{\fd}} \|q_L\|_{L_t^2 \dot{B}_{2,1}^{\fd+1}} +\|q_n\|_{L_t^\infty \dot{B}_{2,1}^{\fd}} \|u_L\|_{L_t^1 \dot{B}_{2,1}^{\fd+1}}\\
+tC(\|q_0\|_{L^\infty}) \|q_n\|_{L_t^\infty \dot{B}_{2,1}^{\fd}} +C\|u_n\|_{L_t^2 \dot{B}_{2,1}^{\fd}} \|u_L\|_{L_t^2 \dot{B}_{2,1}^{\fd}} +C(\|q_0\|_{L^\infty}) \|q_n\|_{L_t^\infty \dot{B}_{2,1}^{\fd}} \|u_L\|_{L_t^1 \dot{B}_{2,1}^{\fd+1}}\\
\leq C_0 \Fs \left( \|q_L\|_{L_t^2 \dot{B}_{2,1}^{\fd+1}} +\|u_L\|_{L_t^1 \dot{B}_{2,1}^{\fd+1}} +\|u_L\|_{L_t^2 \dot{B}_{2,1}^{\fd}}\right).
\end{multline}
Thanks to condition \eqref{condbetaT2}, we obtain the second point.
\begin{rem}
\sl{This is here (precisely for $u_n\cdot \nabla q_L$) that we needed the additional regularity on the initial data.}
\end{rem}
The proof of $(3)$ immediately follows for $T$ small enough so that:
$$
\frac{\nu_0 b_*}{4} \|u_L\|_{L_T^1 \dot{B}_{2,1}^{s+1}} +\gamma^* \|q_L\|_{L^1_T\dot{B}_{\frac{1}{\ee}}^{s+2,s}} \leq \beta.
$$
For the estimates on the time derivative of $q_{n+1}$ (point $(4)$) we simply use the first equation of the iterative scheme and the estimates from $\mathcal{P}(n)$. This implies the result for $J(q_n)$ thanks to composition estimates (see for example \cite{Dbook} or \cite{Dinv} lemma 1.6) as:
$$
\d_t J(q_n)= -\frac{\d_t q_n}{(1+q_n)^2}.
$$
For point $(6)$, we simply write:
$$
\|(I_d-\dot{S}_m)q_{n+1}\|_{\tilde{L}^\infty_T\dot{B}_{2,1}^\fd} \leq \|\overline{q}_{n+1}\|_{\tilde{L}^\infty_T\dot{B}_{2,1}^\fd} +\|(I_d-\dot{S}_m)q_L\|_{\tilde{L}^\infty_T\dot{B}_{2,1}^\fd}.
$$
Using the estimates on $q_L$ and on $\overline{q}_{n+1}$(given by $(2)_{n+1}$, we obtain $(6)$. The following point is simply given by estimating the equation:
$$
\d_t (q_L-q_0)= -\div u_L,
$$
which implies for $T$ small enough
$$
\|q_L-q_0\|_{\tilde{L}_t^\infty \dot{B}_{2,1}^{\fd+1}} \leq \|u_L\|_{L_t^1 \dot{B}_{2,1}^{\fd+1}} \leq \beta,
$$
so that, writing $q_{n+1}-q_0 =\overline{q}_{n+1} +(q_L-q_0)$ we obtain $(7)$. For the last estimates on $J(q_{n+1})$: the first one is obtained thanks to composition estimates in Besov spaces on:
$$
J(q_{n+1}) -1=-\frac{q_{n+1}}{1+q_{n+1}}.
$$
Next we estimate the equation on $J(q_{n+1})-J(q_0)$ (we refer to \cite{Dtruly} proposition 8 for transport estimates):
\begin{multline}
\d_t \left(J(q_{n+1})-J(q_0)\right) +u_n\cdot \nabla \left(J(q_{n+1}) -J(q_0)\right)\\
=\left(\frac{1+q_n}{(1+q_{n+1})^2}-1+1 \right) \div u_{n+1} -u_n\cdot \nabla(J(q_0)-1),
\end{multline}
Thanks to the classical transport estimates, $\mathcal{P}(n)$ and the Gronwall lemma, we obtain that:
$$
\|J(q_{n+1})-J(q_0)\|_{\tilde{L}^\infty_t\dot{B}_{2,1}^\fd} \leq 2 C_0 t^{\frac{1}{2}} e^{2CC_0 t^{\frac{1}{2}}} \left[ 2C_0 (1+2CC_0)^2 +\|J(q_0)-1\|_{\dot{B}_{2,1}^{\fd+1}}\right].
$$
Then for $T$ small enough as stated in the lemma we get the result. The last point is obtained writing that:
$$
(I_d-\dot{S}_m)\left(J(q_{n+1})-1\right) =(I_d-\dot{S}_m)\left(J(q_{n+1})-J(q_0)\right) +(I_d-\dot{S}_m)\left(J(q_0)-1\right),
$$
and using \eqref{condm} and $(8)_n$. $\blacksquare$

The next step in the proof of the existence theorem is then to prove that $(q_n,u_n)$ is a Cauchy sequence in $E_{1/\ee}^{\fd}(T)$ for $T$ small enough as below. For that we define:
$$
(\delta q_n, \delta u_n) =(q_{n+1} -q_n, u_{n+1} -u_n) =(\overline{q}_{n+1} -\overline{q}_n, \overline{u}_{n+1} -\overline{u}_n).
$$
It solves the following system:
$$
\begin{cases}
\begin{aligned}
&\d_t \delta q_n+ u_n\cdot\nabla \delta q_n+ (1+q_n)\div \delta u_n= \delta F,\\
&\d_t \delta u_n+ u_n\cdot\nabla \delta u_n-\frac{1}{1+q_n}\cA \delta u_n-\kappa \frac{\phie*\nabla \delta q_n-\nabla \delta q_n}{\ee^2}= \delta G,\\
&(\delta q_n, \delta u_n)_{|t=0}=(0, 0),
\end{aligned}
\end{cases}
$$
with
$$
\begin{cases}
\delta F= -\delta u_{n-1}.\nabla q_n - \delta q_{n-1} \div u_n,\\
\delta G= -\nabla(H(1+q_n)-H(1+q_{n-1})) -\delta u_{n-1} \cdot u_n +\left(\frac{1}{1+q_n} -\frac{1}{1+q_{n-1}} \right) \cA u_n.
\end{cases}
$$
Then using theorem \ref{thestimbc} we prove that for sufficiently small $\beta>0$ and $T>0$,
$$
\|(\delta q_n, \delta u_n)\|_{E_{1/\ee}^\fd} \leq \beta^{\frac{1}{3}}<1.
$$
This part is classical (see for example \cite{Dtruly, Has7, Has8}) so we shall not give details (we also refer to the next part, the present computations are close to the ones to estimate the convergence rate). This implies that $(q_n,u_n)$ is a Cauchy sequence in $E_{1/\ee}^{\fd}(T)$ for $T$ small enough, and then this gives a solution for system $(NSR_\ee)$in the case of smooth initial data: $q_0 \in \dot{B}_{2,1}^\fd \cap \dot{B}_{2,1}^{\fd+1}$ and $u_0 \in \dot{B}_{2,1}^{\fd-1} \cap \dot{B}_{2,1}^\fd$.

In the case $q_0 \in \dot{B}_{2,1}^\fd$ and $u_0 \in \dot{B}_{2,1}^{\fd-1}$ we do as in \cite{Dtruly}: we begin by regularizing the initial data:
$$
(q_0^p, u_0^p)= (\dot{S}_p q_0, \dot{S}_p u_0),
$$
and thanks to what precedes, we get a sequence of smooth solutions $(q_p, u_p)_{p\in \N}$ of $(NSR_\ee)$ on a bounded intervall $[0,T]$ (uniformly in $p$). We then use classical arguments (estimates on the time derivatives, compactness, Fatou lemma) to prove that the sequence converges, up to an extraction, towards a solution of $(NSR_\ee)$ with initial data $(q_0, u_0)$ and satisfying the energy estimates. We refer to \cite{Dtruly} for details: in our case the method is very close except that we mainly use our a priori estimates from theorem \ref{thestimbc}.
\begin{rem}
\sl{We could simplify the proof of $(1)_{n+1}$ by using also $(1)_n$ for $s=\fd+1$, then it would not have been necessary to ask for the condition on $\beta$ in \eqref{condbetaT}, but we choosed to present it this way because for $(NSK)$ this is the only way to do as we do not ask for more regularity.}
\end{rem}
\begin{rem}
The proof for (NSK) follows the same lines but is easier: thanks to the greater regularity of $q$, we can immediately prove that $(q_n,u_n)$ is a Cauchy sequence. There is no need to regularize, and therefore $\mathcal{P}(n)$ is simpler: (1) is only proved for $s=\fd$.
\end{rem}
\begin{rem}
A very quick method consists in using what is proved for $(NSC)$ to our system $(NSR_\ee)$. The problem is that the method from \cite{Dtruly} only provides a lower bound for the maximal lifespan $T_\ee^*$ in $\mathcal{O}(\ee^2)$, but we can use our a priori estimates to bound $T_\ee^*$ from below by the same small $T>0$ exhibed in what precedes. This method is valid for the non-local systems but fails for $(NSK)$.
\end{rem}

\subsection{Uniqueness}

The method is similar to what is done in \cite{Dtruly}, except that we use our a priori estimates. As $\ee$ is fixed we can also simply use the estimates from \cite{Dtruly} to get the uniqueness.

\section{Proof of the convergence}

Let us now turn to the convergence rate. Once we obtained \eqref{estimgs} the end of the proof is very close to the one from \cite{CH}. Let $(q_0, u_0) \in \dot{B}_{2,1}^\fd \times\dot{B}_{2,1}^{\fd-1}$ and $(\qe,\ue)$ and $(q,u)$ the corresponding solutions of systems $(NSR_\ee)$ and $(NSK)$ both of them existing on $[0,T]$ where $T$ is provided by theorem \ref{thexistR}:
$$
\begin{cases}
\begin{aligned}
&\d_t \qe+ \ue\cdot\nabla \qe+ (1+\qe)\div \ue=0,\\
&\d_t \ue+ \ue\cdot\nabla \ue -\frac{1}{1+\qe}\cA \ue+\nabla (H(1+\qe)-H(1))-\kappa \frac{\phie *\nabla \qe -\nabla \qe}{\ee^2}=0,\\
\end{aligned}
\end{cases}
$$
and let us write $(q,u)$ as a solution of system $(NSR_\ee)$ with an additionnal external force:
$$
\begin{cases}
\begin{aligned}
&\d_t q+ u\cdot\nabla q+ (1+q)\div u=0,\\
&\d_t u+ u\cdot\nabla u -\frac{1}{1+q}\cA u +\nabla (H(1+q)-H(1)) -\kappa \frac{\phie *\nabla q -\nabla q}{\ee^2}=R_\ee,\\
\end{aligned}
\end{cases}
$$
where the remainder is given by $R_\ee\overset{def}{=}-\frac{\kappa}{\ee^2}\nabla(\phi_\ee*q-q-\ee^2 \Delta q)$. Let us define $(\dq, \du)=(\qe-q, \ue-u)$, who solves:
\begin{equation}
\begin{cases}
\begin{aligned}
&\d_t \dq+ \ue \cdot\nabla \dq+ (1+\qe)\div \du = \delta F_1 + \delta F_2,\\
&\d_t \du+ \ue\cdot\nabla \du -\frac{1}{1+\qe}\cA \du -\kappa \frac{\phie *\nabla \dq -\nabla \dq}{\ee^2}=-R_\ee +\delta G_1 +\delta G_2 +\delta G_3,\\
&(\dq, \du)_{|t=0}=(0,0),
\end{aligned}
\end{cases}
\end{equation}
where
$$
\begin{cases}
\delta F_1 =-\du \cdot \nabla q,\\
\delta F_2 =-\dq. \div u,\\
\end{cases}
\mbox{ and }
\begin{cases}
\delta G_1 =-\du \cdot \nabla u,\\
\delta G_2 = \left(\frac{1}{1+\qe} -\frac{1}{1+q} \right) \cA u,\\
\delta G_3 = -\nabla (H(1+\qe)-H(1+q)).
\end{cases}
$$
\begin{rem}
\sl{We emphasize that in the external force terms from the first equation, there is no such term as $q \cdot \div \du$, which was one of the aims of our methods.}
\end{rem}
Applying \eqref{estimgs} to this system for $s=\fd -h$ for $h\in ]0,1]\cap ]0, d-1[$ (we refer to \cite{CH} for explainations of such a choice, mainly due to paraproduct and remainder endpoints) and:
$$
b=\frac{1}{1+\qe}=J(\qe), \quad c=1+\qe, \quad v=\ue,
$$
we obtain
\begin{multline}
g^{\fd-h}(\dq,\du)(t)+\frac{\nu_0 b_*}{4} \|\du\|_{L_t^1 \dot{B}_{2,1}^{\fd-1+1}} +\gamma_* \|\dq\|_{L_t^1 \dot{B}_{\frac{1}{\ee}}^{\fd-h+2, \fd-h}} \leq g^{\fd-h}(0,0)(0)\\
+\Fs \int_0^t \left(\|\delta F_1+\delta F_2\|_{\dot{B}_{2,1}^s} +\|R_\ee +\delta G_1 +\delta G_2 +\delta G_3\|_{\dot{B}_{2,1}^{s-1}}\right) d\tau\\
+ \Fs \int_0^t g^{\fd-h}(\dq,\du)(\tau) \Bigg[ 2^m \Big(\|\d_t b\|_{\dot{B}_{2,1}^{\fd-1}} +\|\d_t c\|_{\dot{B}_{2,1}^{\fd-1}}\Big)\\
+(1+\|b-1\|_{\dot{B}_{2,1}^\fd} +\|c-1\|_{\dot{B}_{2,1}^\fd})^2 \Big(\|\nabla \ue\|_{\dot{B}_{2,1}^\fd} +2^m \|\ue\|_{\dot{B}_{2,1}^\fd} +2^{2m} +\|\ue\|_{\dot{B}_{2,1}^\fd}^2\Big)\Bigg] d\tau\\
+\Fs \int_0^t \left(\|(I_d-\dot{S}_m)(b-1)\|_{\dot{B}_{2,1}^\fd} +\|(I_d-\dot{S}_m)(c-1)\|_{\dot{B}_{2,1}^\fd}\right) \|\ue\|_{\dot{B}_{2,1}^{s+1}} d\tau.
\label{estimgdh}
\end{multline}
Thanks to Theorem \ref{thexistR} we have:
\begin{equation}
\|c-1\|_{\tilde{L}_t^\infty \dot{B}_{2,1}^\fd} =\|\qe\|_{\tilde{L}_t^\infty \dot{B}_{2,1}^\fd} \leq C_0',
\end{equation}
if $T$ is small enough. Thanks to composition estimates (see for example \cite{Dbook}, \cite{Dinv} lemma 1.6 or \cite{CH} appendix), we get:
\begin{equation}
\|b-1\|_{\tilde{L}_t^\infty \dot{B}_{2,1}^\fd} =\|\frac{\qe}{1+\qe}\|_{\tilde{L}_t^\infty \dot{B}_{2,1}^\fd} \leq C(\|\qe\|_{L^\infty}) \|\qe\|_{\tilde{L}_t^\infty \dot{B}_{2,1}^\fd} \leq C_0'.
\end{equation}
Moreover, thanks to the equation on $\qe$,
$$
\|\d_t c\|_{\dot{B}_{2,1}^{\fd-1}} =\|\d_t \qe\|_{\dot{B}_{2,1}^{\fd-1}} \leq (1+\|\qe\|_{\dot{B}_{2,1}^\fd}) \|\ue\|_{\dot{B}_{2,1}^\fd},
$$
and by interpolation we obtain for all $t\leq T$:
\begin{equation}
\|\d_t c\|_{L_t^2 \dot{B}_{2,1}^{\fd-1}} \leq (1+\|\qe\|_{L_t^\infty\dot{B}_{2,1}^\fd}) \|\ue\|_{L_t ^\infty\dot{B}_{2,1}^{\fd-1}}^{\frac{1}{2}} \|\ue\|_{L_t ^1\dot{B}_{2,1}^{\fd+1}}^{\frac{1}{2}} \leq C_0' V_\ee (t)^{\frac{1}{2}},
\end{equation}
with
$$
V_\ee (t) =\int_0^t \|\ue\|_{\dot{B}_{2,1}^{\fd+1}} d\tau.
$$
Similarly, we obtain
\begin{multline}
\|\d_t b\|_{\dot{B}_{2,1}^{\fd-1}} =\|\left[ 1+\Big(\frac{1}{(1+\qe)^2}-1\Big)\right] \d_t \qe \|_{\dot{B}_{2,1}^{\fd-1}}\\
\leq (1+\|\frac{1}{(1+\qe)^2}-1\|_{\dot{B}_{2,1}^\fd}) \|\d_t \qe\|_{\dot{B}_{2,1}^{\fd-1}} \leq (1+C(\|\qe\|_{L^\infty}) \|\qe\|_{\tilde{L}_t^\infty \dot{B}_{2,1}^\fd}) \|\d_t \qe\|_{\dot{B}_{2,1}^{\fd-1}}.
\end{multline}
Collecting these estimates and thanks to the estimates on $(q,u)$ and $(\qe, \ue)$ given by theorem \ref{thexistR}, we conclude that there exists a constant $C_0'>1$ only depending on the initial data such that:
\begin{equation}
\begin{cases}
\|b-1\|_{\tilde{L}_t^\infty \dot{B}_{2,1}^\fd} +\|c-1\|_{\tilde{L}_t^\infty \dot{B}_{2,1}^\fd} \leq C_0',\\
\|\d_t b\|_{L_t^2 \dot{B}_{2,1}^{\fd-1}} +\|\d_t c\|_{L_t^2 \dot{B}_{2,1}^{\fd-1}} \leq C_0' V_\ee(t)^\frac{1}{2}.
\end{cases}
\end{equation}
As $\qe$ satisfies a transport equation, thanks to Proposition 9 from \cite{Dtruly} we can write:
\begin{multline}
\|(I_d-\dot{S}_m)(c-1)\|_{\tilde{L}_t^\infty \dot{B}_{2,1}^\fd} =\|(I_d-\dot{S}_m)\qe\|_{\tilde{L}_t^\infty \dot{B}_{2,1}^\fd}\\
\leq \|(I_d-\dot{S}_m) q_0\|_{\dot{B}_{2,1}^\fd}  +(1+\|q_0\|_{\dot{B}_{2,1}^\fd})(e^{C V_\ee (t)}-1).
\label{estimHFJq}
\end{multline}
Similarly, as $b=J(\qe)$ satisfies:
$$
\d_t J(\qe) +\ue \cdot \nabla J(q)= J(q) \div \ue,
$$
we obtain:
\begin{multline}
\|(I_d-\dot{S}_m)(b-1)\|_{\tilde{L}_t^\infty \dot{B}_{2,1}^\fd} =\|(I_d-\dot{S}_m)(J(\qe)-1)\|_{\tilde{L}_t^\infty \dot{B}_{2,1}^\fd}\\
\leq \|(I_d-\dot{S}_m) (J(q_0)-1)\|_{\dot{B}_{2,1}^\fd}  +(1+\|J(q_0)-1\|_{\dot{B}_{2,1}^\fd})(e^{C V_\ee (t)}-1).
\end{multline}
If $m$ is large enough, and $T$ is small enough (both of them only depending on the initial data and the physical parameters and not on $\ee$), we have:
$$
\|(I_d-\dot{S}_m)(b-1)\|_{\tilde{L}_t^\infty \dot{B}_{2,1}^\fd} +\|(I_d-\dot{S}_m)(b-1)\|_{\tilde{L}_t^\infty \dot{B}_{2,1}^\fd} \leq \frac{\gamma_*}{16 \Fs},
$$
so that \eqref{estimgdh} reduces to ($m$ is now fixed)
\begin{multline}
g^{\fd-h}(\dq,\du)(t)+\frac{ 3\nu_0 b_*}{16} \|\du\|_{L_t^1 \dot{B}_{2,1}^{\fd-1+1}} +\gamma_* \|\dq\|_{L_t^1 \dot{B}_{\frac{1}{\ee}}^{\fd-h+2, \fd-h}}\\
\leq \Fs \int_0^t \left(\|\delta F_1+\delta F_2\|_{\dot{B}_{2,1}^s} +\|R_\ee +\delta G_1 +\delta G_2 +\delta G_3\|_{\dot{B}_{2,1}^{s-1}}\right) d\tau\\
+ \Fs \int_0^t g^{\fd-h}(\dq,\du)(\tau) \Big(\|\nabla \ue\|_{\dot{B}_{2,1}^\fd} +2^m \|\ue\|_{\dot{B}_{2,1}^\fd} +2^{2m} +\|\ue\|_{\dot{B}_{2,1}^\fd}^2\Big)\Bigg] d\tau.\\
\label{estimgdh2}
\end{multline}
The last step is now to estimate the external forcing terms. This is similar to \cite{CH} so we will skip details: using Besov product laws and interpolation we obtain:
\begin{multline}
\|\delta F_1\|_{\dot{B}_{2,1}^{\fd-h}} \leq C\|\du\|_{\dot{B}_{2,1}^{\fd-h}} \|q\|_{\dot{B}_{2,1}^{\fd+1}}\\
\leq \frac{\nu_0 b_*}{16 \Fs} \|\du\|_{\dot{B}_{2,1}^{\fd-h+1}} +\Fs \|\du\|_{\dot{B}_{2,1}^{\fd-h-1}} \|q\|_{\dot{B}_{2,1}^\fd} \|q\|_{\dot{B}_{2,1}^{\fd+2}},
\end{multline}
and
$$
\begin{cases}
\|\delta F_2\|_{\dot{B}_{2,1}^{\fd-h}} \leq C \|\dq\|_{\dot{B}_{2,1}^{\fd-h}} \|\nabla u\|_{\dot{B}_{2,1}^\fd},\\
\|\delta G_1\|_{\dot{B}_{2,1}^{\fd-h-1}} \leq C \|\du\|_{\dot{B}_{2,1}^{\fd-h-1}} \|\nabla u\|_{\dot{B}_{2,1}^\fd},\\
\|\delta G_2\|_{\dot{B}_{2,1}^{\fd-h-1}} \leq C \|I(\qe)-I(q)\|_{\dot{B}_{2,1}^{\fd-h}} \|\cA u\|_{\dot{B}_{2,1}^\fd} \leq C_0' \|\dq\|_{\dot{B}_{2,1}^{\fd-h}} \|\cA u\|_{\dot{B}_{2,1}^\fd},\\
\|\delta G_3\|_{\dot{B}_{2,1}^{\fd-h-1}} \leq C(\|\qe\|_{L^\infty}, \|q\|_{L^\infty}) \|\dq\|_{\dot{B}_{2,1}^{\fd-h}}.
\end{cases}
$$
Plugging this into \eqref{estimgdh2} and using the Gronwall lemma allows us to write for all $t\leq T$:
\begin{multline}
g^{\fd-h}(\dq,\du)(t)+\frac{ 3\nu_0 b_*}{16} \|\du\|_{L_t^1 \dot{B}_{2,1}^{\fd-h+1}} +\gamma_* \|\dq\|_{L_t^1 \dot{B}_{1/\ee}^{\fd-h+2, \fd-h}}\\
\leq \Fs \|R_\ee\|_{L_t^1 \dot{B}_{2,1}^{s-1}} e^{\Fs \int_0^t \Big(\|\nabla \ue\|_{\dot{B}_{2,1}^\fd} +2^m \|\ue\|_{\dot{B}_{2,1}^\fd} +1 +2^{2m} +\|\ue\|_{\dot{B}_{2,1}^\fd}^2 +\|\nabla u\|_{\dot{B}_{2,1}^\fd} +\|q\|_{\dot{B}_{2,1}^\fd} \|q\|_{\dot{B}_{2,1}^{\fd+2}} \Big) d\tau}\\
\end{multline}
We can estimate $R_\ee$ thanks to Corollary 1 from \cite{CH} (section 4): there exists a constant $C_h>0$ such that:
$$
\|R_\ee\|_{L_T^1 \dot{B}_{2,1}^{\fd-h-1}} \leq C_h \kappa \ee^h \|q\|_{L_T^1 \dot{B}_{2,1}^{\fd+2}}.
$$
We conclude using the estimates on $q$ given by theorem \ref{thexistR}.
\begin{rem}
\sl{As in \cite{CH, Corder}, in the case $h=0$, using the same estimates we have:
$$
\|R_\ee\|_{L_t^1 \dot{B}_{2,1}^{\fd-1}}\leq \kappa \sum_{j\in \Z} 2^{j\fd}\frac{e^{-\ee^2 C^2 2^{2j}}-1+\ee^2 C^2 2^{2j}}{\ee^2} \|\ddj q\|_{L_t^1 L^2},
$$
and thanks to the Lebesgue theorem it goes to zero as $\ee$ goes to zero.
}
\end{rem}

\section{Variable coefficients}

In this section we explain how our results can be extended for density dependant viscosity coefficients. Let us recal that in this case the diffusion operator writes:
$$
\cA u = \div (2\mu(\rho) Du) +\nabla (\lambda(\rho) \div u),
$$
where $2Du =^t\nabla u +\nabla u$. Then we obtain:
$$
\frac{1}{\rho} \cA u =\frac{\mu(\rho)}{\rho} \Delta u +\frac{\mu(\rho)+\lambda(\rho)}{\rho} \nabla \div u + 2Du\cdot \frac{\mu'(\rho)}{\rho}
 \nabla \rho +\div u \cdot \frac{\lambda'(\rho)}{\rho} \nabla \rho.$$
The last two terms are of the form $\nabla u \cdot \nabla L(\rho)$ and can be dealt the same way as we did in \cite{CHSW} writing:
$$
\nabla u \cdot \nabla L(1+q)= \nabla u \cdot \nabla \left(L(1+q)-L(1+\dot{S}_m q) \right) +\nabla u \cdot \nabla \left(L(1+\dot{S}_m q) -L(1)\right).
$$
For the first two terms, we have the diffusion operator:
$$
\mu_1(1+q) \Delta u +\Big(\mu_1(1+q)+\lambda_1(1+q)\Big) \nabla \div u,
$$
with $\mu_1(x)=\mu(x)/x$, $\lambda_1(x)=\lambda(x)/x$. Instead of dealing with the operator $b(1+q) (\mu \Delta u +(\lambda+\mu)\nabla \div u)$ in the proof of theorem \ref{thestimbc}, we have to deal with $b_1(1+q) \Delta u +b_2(1+q) \nabla \div u$ where we will require estimates on $\d_t b_1$, $\d_t b_2$, $\d_t c$, $c-1$, $b_1 -b_1(1)$ and $b_2-b_2(1)$.

The last point we need to carefully check is the transport estimates involving $\mu_1(1+q)$ or $\lambda_1(1+q)$ like in (5) (8) (9) (10) from $\mathcal{P}(n)$ or \eqref{estimHFJq}. This relies on the composition estimates for Besov norms (see lemma 1.6 from \cite{Dinv} or \cite{Dbook}) which require that
$$
\mu(x)/x, \quad \mu'(x), \quad \lambda(x)/x, \quad \lambda'(x) \in W_{local}^{[\fd]+2, \infty}.
$$
We emphasize that for special cases $\mu(\rho)= \mu \rho$ and $\lambda(\rho)= \lambda \rho$ we only need the estimates from \cite{CH}.

\section{Other Besov settings}

Our method relies on innerproducts based on $L^2$ and cannot be adapted in the setting of $\dot{B}_{p,r}^s$ spaces with $p\neq 2$. In the case of particular coefficients (viscosity, capillarity) we can use the decoupled methods (we refer to \cite{Dtruly, CHSW, Has7, Has8}) and obtain results in this case.

For a summation index $r\neq 1$, the methods we present here can be adapted in the same way as in \cite{CHSW, Has7, Has8}.

\section{Appendix}

\subsection{Besov spaces}

\subsubsection{Littlewood-Paley theory}

The Fourier transform of $u$ with respect to the space variable will be denoted by $\mathcal{F}(u)$ or $\hat{u}$. 
In this subsection we will briefly state (as in \cite{CD}) classical definitions and properties concerning the homogeneous dyadic decomposition with respect to the Fourier variable. We will recall some classical results and we refer to \cite{Dbook} (Chapter 2) for proofs and a complete presentation of the theory.

To build the Littlewood-Paley decomposition, we need to fix a smooth radial function $\chi$ supported in (for example) the ball $B(0,\frac{4}{3})$, equal to 1 in a neighborhood of $B(0,\frac{3}{4})$ and such that $r\mapsto \chi(r.e_r)$ is nonincreasing over $\R_+$. So that if we define $\varphi(\xi)=\chi(\xi/2)-\chi(\xi)$, then $\varphi$ is compactly supported in the annulus $\{\xi\in \R^d, c_0=\frac{3}{4}\leq |\xi|\leq C_0=\frac{8}{3}\}$ and we have that,
\begin{equation}
 \forall \xi\in \R^d\setminus\{0\}, \quad \sum_{l\in\Z} \varphi(2^{-l}\xi)=1.
\label{LPxi}
\end{equation}
Then we can define the \textit{dyadic blocks} $(\ddl)_{l\in \Z}$ by $\ddl:= \varphi(2^{-l}D)$ (that is $\hat{\ddl u}=\varphi(2^{-l}\xi)\hat{u}(\xi)$) so that, formally, we have
\begin{equation}
u=\Sum_l \ddl u
\label{LPsomme} 
\end{equation}
We can now define the homogeneous Besov spaces used in this article:
\begin{defi}
\label{LPbesov}
 For $s\in\R$ and  
$1\leq p,r\leq\infty,$ we set
$$
\|u\|_{\dot B^s_{p,r}}:=\bigg(\sum_{l} 2^{rls}
\|\ddl u\|^r_{L^p}\bigg)^{\frac{1}{r}}\ \text{ if }\ r<\infty
\quad\text{and}\quad
\|u\|_{\dot B^s_{p,\infty}}:=\sup_{l} 2^{ls}
\|\ddl u\|_{L^p}.
$$
We then define the space $\dot B^s_{p,r}$ as the subset of  distributions $u\in {\cS}'_h$ such that $\|u\|_{\dot B^s_{p,r}}$ is finite.
\end{defi}
Once more, we refer to \cite{Dbook} (Chapter $2$) for properties of the inhomogeneous and homogeneous Besov spaces.

In this paper, we mainly work with functions or distributions depending on both the time variable $t$ and the space variable $x$, and we introduce the spaces $L^p([0,T];X)$ (resp. $\cC([0,T];X)$).

The Littlewood-Paley decomposition enables us to work with spectrally localized (hence smooth) functions rather than with rough objects. We naturally obtain bounds for each dyadic block in spaces of type $L^\rho_T L^p.$  Going from those type of bounds to estimates in  $L^\rho_T \dot B^s_{p,r}$ requires to perform a summation in $\ell^r(\Z).$ When doing so however, we \emph{do not} bound the $L^\rho_T \dot B^s_{p,r}$ norm for the time integration has been performed \emph{before} the $\ell^r$ summation.
This leads to the following notation (after J.-Y. Chemin and N. Lerner in \cite{CL}):

\begin{defi}\label{d:espacestilde}
For $T>0,$ $s\in\R$ and  $1\leq r,\rho\leq\infty,$
 we set
$$
\|u\|_{\tilde L_T^\rho \dot B^s_{p,r}}:=
\bigl\Vert2^{js}\|\ddq u\|_{L_T^\rho L^p}\bigr\Vert_{\ell^r(\Z)}.
$$
\end{defi}
\medbreak
Let us now recall a few nonlinear estimates in Besov spaces. Formally, any product of two distributions $u$ and $v$ may be decomposed into 
\begin{equation}\label{eq:bony}
uv=T_uv+T_vu+R(u,v), \mbox{ where}
\end{equation}
$$
T_uv:=\sum_l\dot S_{l-1}u\ddl v,\quad
T_vu:=\sum_l \dot S_{l-1}v\ddl u\ \hbox{ and }\ 
R(u,v):=\sum_l\sum_{|l'-l|\leq1}\ddl u\,\dot\Delta_{l'}v,
$$
where for $j\in \Z$, the operator $\dot{S}_j$ is defined by:
$$\dot{S}_j u =\Sum_{l\leq j-1} \ddl u = \chi(2^{-j}D) u.$$
The above operator $T$ is called a ``paraproduct'' whereas $R$ is called a "remainder''. The decomposition \eqref{eq:bony} has been introduced by J.-M. Bony in \cite{BJM}. We refer to \cite{Dbook} for properties, and also to \cite{CH} or \cite{CHVP} for paraproduct and remainder estimates for external force terms.

In this article we will frequently use the following estimates (we refer to \cite{Dbook} Section 2.6, \cite{Dinv}, \cite{Has1} for general statements, more properties of continuity for the paraproduct and remainder operators, sometimes adapted to $\tilde L_T^\rho \dot B^s_{p,r}$ spaces): under the same assumptions there exists a constant $C>0$ such that:
\begin{equation}
\|\dot{T}_u v\|_{\dot{B}_{2,1}^s}\leq C \|u\|_{L^\infty} \|v\|_{\dot{B}_{2,1}^s}\leq C \|u\|_{\dot{B}_{2,1}^\fd} \|v\|_{\dot{B}_{2,1}^s},
\label{estimbesov}
\end{equation}
$$
\|\dot{T}_u v\|_{\dot{B}_{2,1}^{s+t}}\leq C\|u\|_{\dot{B}_{\infty,\infty}^t} \|v\|_{\dot{B}_{2,1}^s} \leq C\|u\|_{\dot{B}_{2,1}^{t+\fd}} \|v\|_{\dot{B}_{2,1}^s} \quad (t<0),
$$
$$\|\dot{R}(u,v)\|_{\dot{B}_{2,1}^{s_1+s_2}} \leq C\|u\|_{\dot{B}_{\infty,\infty}^{s_1}} \|v\|_{\dot{B}_{2,1}^{s_2}} \leq C\|u\|_{\dot{B}_{2,1}^{s_1+\fd}} \|v\|_{\dot{B}_{2,1}^{s_2}} \quad (s_1+s_2>0),
$$
$$
 \|\dot{R}(u,v)\|_{\dot{B}_{2,1}^{s_1+s_2-\fd}} \leq C\|\dot{R}(u,v)\|_{\dot{B}_{1,1}^{s_1+s_2}} \leq C\|u\|_{\dot{B}_{2,1}^{s_1}} \|v\|_{\dot{B}_{2,1}^{s_2}} \quad (s_1+s_2>0).
$$
Let us now turn to the composition estimates. We refer for example to \cite{Dbook} (Theorem $2.59$, corollary $2.63$)):
\begin{prop}
\sl{\begin{enumerate}
 \item Let $s>0$, $u\in \dot{B}_{2,1}^s\cap L^{\infty}$ and $F\in W_{loc}^{[s]+2, \infty}(\R^d)$ such that $F(0)=0$. Then $F(u)\in \dot{B}_{2,1}^s$ and there exists a function of one variable $C_0$ only depending on $s$, $d$ and $F$ such that
$$
\|F(u)\|_{\dot{B}_{2,1}^s}\leq C_0(\|u\|_{L^\infty})\|u\|_{\dot{B}_{2,1}^s}.
$$
\item If $u$ and $v\in\dot{B}_{2,1}^\fd$ and if $v-u\in \dot{B}_{2,1}^s$ for $s\in]-\fd, \fd]$ and $F\in W_{loc}^{[s]+3, \infty}(\R^d)$, then $F(v)-F(u)$ belongs to $\dot{B}_{2,1}^s$ and there exists a function of two variables $C$ only depending on $s$, $d$ and $G$ such that
$$
\|F(v)-F(u)\|_{\dot{B}_{2,1}^s}\leq C(\|u\|_{L^\infty}, \|v\|_{L^\infty})\left(|G'(0)| +\|u\|_{\dot{B}_{2,1}^\fd} +\|v\|_{\dot{B}_{2,1}^\fd}\right) \|v-u\|_{\dot{B}_{2,1}^s}.
$$
\end{enumerate}}
\label{estimcompo}
\end{prop}

We end this subsection with two commutator estimates proven in \cite{Dtruly} (we also refer to \cite{Dbook}).

\begin{lem}
\sl{let $\sigma \in ]-\fd, \fd+1]$. There exists a sequence $(c_j)_{j\in \Z}\in l^1(\Z)$ with summation 1, and a constant $C$ only depending on $d$ and $\sigma$ such that for all $j\in \Z$,
$$
\|[v\cdot \nabla, \ddj]a\|_{L^2} \leq C c_j 2^{-j \sigma} \|\nabla v\|_{\dot{B}_{2,1}^\fd} \|a\|_{\dot{B}_{2,1}^\sigma}.
$$
}
\end{lem}

\begin{lem}
\sl{Let $\sigma \in ]-\fd, \fd+1]$ and $h \in ]1-\fd, 1]$, $k\in \{1,...,d\}$ and $\mathcal{R}_j=\ddj (a \d_k w)- \d_k (a\ddj w)$. There exists a constant $C$ only depending on $\sigma, h,d$ such that:
$$
\Sum_{j\in \Z} 2^{j\sigma} \|\mathcal{R}_j\|_{L^2} \leq C \|a\|_{\dot{B}_{2,1}^{\fd+h}} \|w\|_{\dot{B}_{2,1}^{\sigma+1-h}}.
$$
}
\label{estimtruly}
\end{lem}

\subsubsection{Complements for hybrid Besov spaces}

As explained, in the study of the compressible Navier-Stokes system with data in critical spaces, the density fluctuation has two distinct behaviours in low and high frequencies, separated by a frequency threshold. This leads to the notion of hybrid Besov spaces and we refer to R. Danchin in \cite{Dinv} or \cite{Dbook} for general hybrid spaces . In this paper we only will use the following hybrid norms:
\begin{multline}
 \|f\|_{\dot{B}_\beta^{s+2,s}} \overset{def}{=} \sum_{j\in \Z} \min(\beta^2, 2^{2j}) 2^{js} \|\ddj f\|_{L^2}\\
 = \Sum_{j\leq \log_2 \beta} 2^{j(s+2)} \|\ddj f\|_{L^2} +\Sum_{j> \log_2 \beta} \beta^2 2^{js} \|\ddj f\|_{L^2}.
\end{multline}
In this formulation, we obviously remark the threshold frequency $\log_2 \beta$ which separates low (parabolically regularized) and high (damped with coefficient $\beta^2$) frequencies. But as we prove in \cite{CHVP, Corder}, the frequency transition is in fact continuous and the following equivalent formulations show that these norms are completely tailored to our capillary term. We refer to \cite{CHVP} for the case of the first non-local model:
$$
 \|f\|_{\dot{B}_{1/\ee}^{s+2,s}} \sim \sum_{j\in \Z} \frac{1-e^{-c\ee^2 2^{2j}}}{\ee^2}  2^{js} \|\ddj f\|_{L^2}  \sim \|\frac{\phie*f-f}{\ee^2}\|_{\dot{B}_{2,1}^s} 
$$
and we refer to \cite{Corder} for the order parameter model:
$$
\|f\|_{\dot{B}_\aa^{s+2,s}} \sim \sum_{j\in \Z} \frac{2^{2j}}{1+\frac{2^{2j}}{\aa^2}}  2^{js} \|\ddj f\|_{L^2} \sim \|\aa^2 (\psia*f-f)\|_{\dot{B}_{2,1}^s}.
 $$


\end{document}